\documentclass[a4paper]{article}

\usepackage{geometry}

\usepackage{graphicx}
\usepackage{xspace}

\usepackage{hyperref}

\usepackage{tabu}
\usepackage{multirow}
\usepackage{booktabs}
\usepackage{sidecap}
\usepackage[binary-units]{siunitx}

\usepackage{float}

\usepackage{latexsym}
\usepackage{amsmath}
\usepackage{amsfonts}
\usepackage{amssymb}
\usepackage{stmaryrd}
\usepackage{graphicx,epsfig}

\usepackage[thref,amsmath]{ntheorem}

\newtheorem{theorem}{Theorem}
\newtheorem{lemma}{Lemma}
\newtheorem{corollary}{Corollary}
\newtheorem{definition}{Definition}
\newtheorem{fact}{Fact}
\newtheorem{proposition}{Proposition}
\newtheorem{example}{Example}
\newtheorem{observation}{Observation}

\newtheorem{remark}{Remark}

\theoremstyle{nonumberplain}
\theoremheaderfont{\itshape}
\theorembodyfont{\normalfont}
\theoremseparator{.}
\theoremsymbol{}
\newtheorem{proof}{Proof}

\providecommand{\keywords}[1]{\textbf{\textit{Keywords ---}} #1}

\usepackage{authblk}

\newbox{\myorcidaffilbox}
\sbox{\myorcidaffilbox}{\large\includegraphics[height=1.7ex]{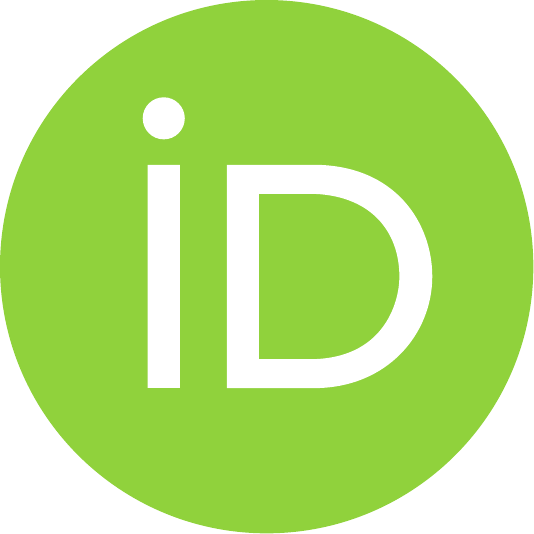}}
\newcommand{\orcid}[1]{%
	\href{https://orcid.org/#1}{\usebox{\myorcidaffilbox}}}

\usepackage{epstopdf}

\usepackage[ddmmyyyy]{datetime}

\usepackage[ruled,noend,noline]{algorithm2e}
\usepackage{yfonts}
\usepackage{url}

\newcommand{\qed}{\hfill \ensuremath{\Box}}

\renewcommand{\epsilon}{\varepsilon}
\renewcommand{\phi}{\varphi}

\def\IN{{\mathbb N}}

\def\EL{{\cal L}}
\def\L{{\cal L}} 
\def\pref{\mathop{{\rm pref}}}

\def\lex{{\rm lex}}

\newcommand{\rank}{{\textit{rank}}}

\newcommand{\pv}{pv}

\newcommand{\op}{\ensuremath{\mathop{\rm op}}}
\newcommand{\flipext}{\ensuremath{\mathop{\rm flipext}}}

\newcommand{\fin}{\text{\rm fin}}

\renewcommand{\epsilon}{\varepsilon}

\newcommand{\PNF}{\mathrm{PNF}}

\newcommand{\tm}{{\textrm{TM}}}

\newcommand{\Fac}{\textit{Fct}}
\newcommand{\Fct}{\textit{Fct}}

\newcommand{\lazyflipext}{\ensuremath{\mathop{\rm lazy\text{-}flipext}}}
\newcommand{\lazyaflipext}{\ensuremath{\mathop{\rm lazy\text{-}flipext}}}

\usepackage{color}

\usepackage{multirow}
\usepackage{tabularx}
\usepackage{booktabs}
\newcolumntype{Y}{>{\raggedleft\arraybackslash}X}

\begin{document}

\title{On Infinite Prefix Normal Words}

\author[1]{Ferdinando Cicalese\orcid{0000-0003-1652-0599}}
\author[1]{Zsuzsanna Lipt\'{a}k\orcid{0000-0002-3233-0691}\thanks{Corresponding author}}
\author[2]{Massimiliano Rossi\orcid{0000-0002-3012-1394}}

\affil[1]{Department of Computer Science, University of Verona, Verona, Italy}
\affil[2]{Department of Computer and Information Science and Engineering, University of Florida, Gainesville, FL, United States}
{
    \makeatletter
	\renewcommand\AB@affilsepx{, \protect\Affilfont}
	\makeatother
    \affil[1]{\{ferdinando.cicalese, zsuzsanna.liptak\}@univr.it}
    \affil[2]{rossi.m@ufl.edu}
}

\date{\bigskip \bigskip {\em Article published in Theoretical Computer Science (2021)} \\
doi: 10.1016/j.tcs.2021.01.015}

\maketitle

\begin{abstract}
Prefix normal words are binary words with the property that no factor has more $1$s than the prefix of the same length. Finite prefix normal words were introduced in [Fici and Lipt\'ak, DLT 2011]. In this paper, we study infinite prefix normal words and explore their relationship to some known classes of infinite binary words. In particular, we establish a connection between prefix normal words and Sturmian words, between prefix normal words and abelian complexity, and between prefix normality and lexicographic order.\footnote{This is an extended version of our paper presented at  SOFSEM 2019~\cite{CLR19}.}
\end{abstract}



\section{Introduction}\label{sec:introduction}

Prefix normal words are binary words where no factor has more $1$s than the prefix of the same length. As an example, the word $11100110101$ is prefix normal, while $11100110110$ is not, since it has a factor of length 5 with four 1s, while the prefix of length 5 has only three 1s. Finite prefix normal words were introduced in~\cite{FL11} and further studied in~\cite{BFLRS_CPM14,BurcsiFLRS17,SWW17,CLR18,BG19,FKNP20,BFLRS20}. 

One motivation for studying prefix normal words comes from the problem of {\em Indexed Binary Jumbled Pattern Matching}~\cite{IJFCS12,BurcsiCFL12_ToCS,MoosaR_JDA12,GiaGrab_IPL13,AmirCLL_ICAPL14,GagieHLW15,ChanL15,CunhaDGWKS17,ADKN20}:  Given a finite word $s$ of length $n$, construct an index in such a way that the following type of queries can be answered efficiently: for two integers $x,y\geq 0$, does $s$ have a factor with $x$ $1$s and $y$ $0$s? As shown in~\cite{FL11,BurcsiFLRS17}, prefix normal words can be used for constructing such an index, via so-called {\em prefix normal forms}.

Prefix normal words have also been shown to form bubble languages~\cite{RSW12,SW12,BFLRS_CPM14}, a family of binary languages with efficiently generable combinatorial Gray codes; the language of prefix normal words has  connections to the Binary Reflected Gray Code~\cite{SWW17}; 
and, recently, prefix normal words also appeared in a graph theoretic  context~\cite{BM18}. 
Indeed, three sequences related to prefix normal words are present in the On-Line Encyclopedia of Integer Sequences (OEIS~\cite{oeis}): A194850 (the number of prefix normal words of length $n$), A238109 (a list of prefix normal words over the alphabet $\{1,2\}$), and A238110 (maximal equivalence class sizes of words with the same prefix normal form).

\medskip

In~\cite{CLR18}, we introduced infinite prefix normal words and analyzed a particular procedure that, given a finite prefix normal word, extends it while preserving the prefix normality property. We showed that the resulting infinite word is ultimately periodic. In this paper, we present a more comprehensive study of infinite prefix normal words, covering several classes of known and well studied infinite words. We will now give a quick tour of the paper (for precise definitions, see Section~\ref{sec:basics}).

\subsection{Our results}

One way of obtaining infinite prefix normal words is by extending finite prefix normal words. We specify two such operations which, in the limit, produce prefix normal words that are extremal with respect to density (Theorem~\ref{thm:extremal}). 

There exist {\em periodic, ultimately periodic,} and {\em aperiodic} infinite prefix normal words: for example, the periodic words $0^{\omega}, 1^{\omega},$ and $(10)^{\omega}$ are prefix normal; the ultimately periodic word $1(10)^{\omega}$ is prefix normal; and so is the aperiodic word
$10100100010000\cdots = \lim_{n\to \infty} 1010^2\cdots 10^{n}.$
The best studied class of aperiodic words are {\em Sturmian words}. 
We show that a Sturmian word $w$ is prefix normal if and only if $w=1c_{\alpha}$ for some $\alpha$, where $c_{\alpha}$ is the characteristic word of slope $\alpha$ (Theorem~\ref{thm:characteristic}).

We show further that every Sturmian word $w$ can be turned into a prefix normal word by prepending a fixed number of $1$s, which only depends on the slope of $w$. This follows from a more general result regarding $c$-balanced words (Lemma~\ref{lemma:c-bal}). For example, the Fibonacci word

\[ \textswab{f} = 0100101001001010010100100101001001\cdots \]

\noindent is not prefix normal, but  the word $1\textswab{f}$ is. 
Two other well-studied aperiodic words are the Thue-Morse word and the Champernowne word. 
The Thue-Morse word

\[ \textswab{t} = 01101001100101101001011001101001\cdots \]

\noindent is not prefix normal but it can be turned into a prefix normal word by prepending two $1$s: $11\textswab{t}$ is prefix normal. On the other hand, the binary Champernowne word

\[ \textswab{c} = 0110111001011101111000100110101011\cdots \]

\noindent which is constructed by concatenating the binary expansions of the integers in ascending order, is not prefix normal and cannot be turned into a prefix normal word by prepending a finite number of $1$s.

We also show that the notion of {\em prefix normal  forms} from~\cite{FL11,BurcsiFLRS17} can be extended to infinite words. These can be used, similarly to the finite case, to encode the {\em abelian complexity} of the original word.
The study of abelian complexity of infinite words was initiated in~\cite{RSZ11}, and continued e.g.\ in~\cite{MadillR13,BFR14,Turek15,CK16,KaboreK17}. We establish a close relationship between the abelian complexity and the prefix normal forms of $w$ (Theorem~\ref{thm:pnfs_abcomp}). We demonstrate how this close connection can be used to derive results about the prefix normal forms of a word $w$. In some cases, such as for Sturmian words and words which are morphic images under the Thue-Morse morphism, we are able to explicitly give the prefix normal forms of the word (Corollary \ref{coro:tm} and Theorem~\ref{thm:sturmian_pnf}).
Conversely, knowing its prefix normal forms allows us to derive results about the abelian complexity of a word.
We also show how to compute the prefix normal forms of words that are {\it binary uniform morphisms}, based on an algorithm from 
\cite{BSSW16} for computing their abelian complexity.

Another class of well-known binary words are {\em Lyndon words}. Notice that the prefix normal condition is different from the Lyndon condition\footnote{For ease of presentation, we are  using Lyndon to mean lexicographically {\em greatest} among its conjugates; this is equivalent to the usual definition up to renaming characters.}:
for finite words, there are words which are both Lyndon and prefix normal (e.g. $110010$), words which are Lyndon but not prefix normal ($11100110110$), words which are prefix normal but not Lyndon ($110101$), and words which are neither ($101100$). We study infinite prefix normal words and their prefix normal forms in the context of lexicographic orderings, and compare them to infinite Lyndon words~\cite{SMDS94} and the max- and min-words of~\cite{Pirillo05} (Corollary~\ref{coro:pirillo}). 

Finally, we give conditions for periodicity and ultimate periodicity of prefix normal words in terms of their minimum density, a parameter  introduced in~\cite{CLR18} (Theorem~\ref{thm:density}).

\subsection{Overview of paper}
The paper is organized as follows. In Section~\ref{sec:basics}, we introduce our terminology and give some simple facts about prefix normal words.
In Section~\ref{sec:properties_of_infinite_pnw}, we compare different operations that generate
 infinite prefix normal words by extending finite prefix normal words.
In Section~\ref{sec:sturmian}, we study the relationship between Sturmian words and prefix normal words. 
Section~\ref{sec:pnf} deals with the connection between prefix normality and 
abelian complexity, and Section~\ref{sec:lex} focuses on the relationship with lexicographic order.
Finally, in Section~\ref{sec:minimum_density}, we analyze the relationship between periodicity and minimum density of prefix normal words.


\section{Basics}\label{sec:basics}

In our definitions and notations, we follow mostly~\cite{Lothaire3}.
A finite (resp.\ infinite) binary word $w$ is a finite (resp.\ infinite) sequence of elements from $\{0,1\}$. Thus an infinite word is a mapping $w:\IN \to \{0,1\}$, where $\IN$ denotes the set of positive integers. We denote the $i$th character of $w$ by $w_i$. Note that we index words starting from $1$. If $w$ is finite, then its length is denoted by $|w|$. The empty word, denoted $\epsilon$, is the unique word of length $0$. The set of binary words of length $n$ is denoted by $\{0,1\}^n$, the set of all finite words by $\{0,1\}^* = \cup_{n\geq 0} \{0,1\}^n$, and the set of infinite binary words by $\{0,1\}^{\omega}$.
For a finite word $u = u_1 \cdots u_n$, we write $u^{\rm rev} = u_n \cdots u_1$ for the reverse of $u$, and for a finite or infinite word $u$, $\overline{u} = \overline{u}_1 \overline{u}_2\cdots$ for the complement of $u$, where $\overline{a} = 1-a$ for $a \in \{0,1\}$.

For two words $u,v$, where $u$ is finite and $v$ is finite or infinite, we write $uv$ for their concatenation. If $w=uxv$, then $u$ is called a prefix, $x$ a factor (or substring), and $v$ a suffix of $w$.
We denote the set of factors of $w$ by $\Fac(w)$ and
its prefix of length $i$ by $\pref_w(i)$, where $\pref_w(0) = \epsilon$. For a finite word $u$, we write $|u|_1$ for the number of $1$s, and $|u|_0$ for the number of $0$s in $u$, and refer to $|u|_1$ as the {\em weight} of $u$. The {\em Parikh vector} of $u$ is $\pv(u)=(|u|_0,|u|_1)$. A word $w$ is called {\em balanced} if for all $u,v\in \Fac(w)$, $|u|=|v|$ implies $||u|_1 - |v|_1| \leq 1$, and {\em $c$-balanced} if $|u|=|v|$ implies $||u|_1 - |v|_1| \leq c$.

For an integer $k\geq 1$ and $u\in \{0,1\}^n$, $u^k$ denotes the $kn$-length word $uuu\cdots u$ ($k$-fold concatenation of $u$) and $u^{\omega}$ the infinite word $uuu\cdots$. An infinite word $w$ is called {\em periodic} if $w = u^{\omega}$ for some non-empty word $u$, and {\em ultimately periodic} if it can be written as $w=vu^{\omega}$ for some $v$ and non-empty $u$.
A word that is neither periodic nor ultimately periodic is called {\em aperiodic}.
We set $0<1$ and denote by $\leq_{\lex}$ the {\em lexicographic order} between words, i.e.\ $u\leq_{\lex} v$ if $u$ is a prefix of $v$ or there is an index $i\geq 1$ s.t.\ $\pref_u(i-1) = \pref_v(i-1)$ and $u_i < v_i$.

For an operation $\op : \{0,1\}^\ast \rightarrow \{0,1\}^\ast$, we denote by $\op^{(i)}$ the $i$th iteration of $\op$. Further, let $\op^\ast(w) = \{\op^{(i)}(w) \mid i \geq 1 \}$ and $\op^{\omega}(w) = \lim_{i \rightarrow \infty}\op^{(i)}(w)$, if it exists.

A {\em binary morphism} $\mu$ is a function $\mu: \{0,1\}^* \to \{0,1\}^*$ such that for all $u,v \in \{0,1\}^*$, $\mu(uv) = \mu(u)\mu(v)$. A binary morphism $\mu$ is called {\em uniform} if $|\mu(0)| = |\mu(1)|$. A {\em fix point} of a morphism $\mu$ is an infinite word $v$ such that $v = \mu^{\omega}(a)$ for some $a\in \{0,1\}$. 

\begin{definition}
Let $w$ be a (finite or infinite) binary word. We define the following functions: 

\begin{itemize}
\item $P_w(i) = |\pref_w(i)|_1$, the {\em weight} of the prefix of length $i$, 
\item $D_w(i) = P_w(i)/i$, the {\em density} of the prefix of length $i$, 
\item $F^1_w(i) = \max \{ |u|_1 : u \in \Fct(w), |u|=i \}$ the maximum number of $1$s in a factor of length $i$, 
\item $f^1_w(i) = \min \{ |u|_1 : u \in \Fct(w), |u|=i \}$, the minimum number of $1$s in a factor of length $i$,
\item $F^0_w(i) = \max \{ |u|_0 : u \in \Fct(w), |u|=i \}$, the maximum number of $0$s in a factor of length $i$, 
\item $f^0_w(i) = \min \{ |u|_0 : u \in \Fct(w), |u|=i \}$, the minimum number of $0$s in a factor of length $i$.
\end{itemize}

\end{definition}

Note that in the context of succinct indexing, the function $P_w(i)$ is often called $\rank_1(w,i)$. We are now ready to define prefix normal words.

\begin{definition}[Prefix normal words]
A (infinite or finite) binary word $w$ is called {\em $1$-prefix normal}, or simply {\em prefix normal}, if $P_w(i) = F^1_w(i)$ for all $i\geq 1$ (for all $1\leq i \leq |w|$ if $w$ is finite). It is called {\em $0$-prefix normal} if $i-P_w(i) = F^0_w(i)$ for all $i\geq 1$ (for all $1\leq i \leq |w|$ if $w$ is finite).
We denote the set of all finite $1$-prefix normal words by $\EL_{\fin}$, the set of all infinite $1$-prefix normal words by $\EL_{\inf}$, and $\EL = \EL_{\fin} \cup \EL_{\inf}$.
\end{definition}

In other words, a word is prefix normal if no factor has more $1$s than the prefix of the same length. Given a binary word $w$, we say that a factor $u$ of $w$ {\em satisfies the prefix normal condition} if $|u|_1 \leq P_w(|u|)$. 

\begin{example} The word $110100110110$ is not prefix normal since the factor $11011$ has four $1$s, which is more than in the prefix $11010$ of length 5. The word $110100110010$, on the other hand, is prefix normal.
The infinite word $(11001)^\omega$ is not prefix normal, because it has $111$ as a factor, which has more $1$s than the prefix of length 3, but the word $(11010)^\omega$ is.
\end{example}

The following facts about infinite prefix normal words are immediate. 

\begin{lemma}\label{lemma:simple pn facts}
\begin{enumerate}
\item For all $u \in \EL_{\fin}$, the word $w = u0^{\omega} \in \EL_{\inf}$.
\item Let $w\in \{0,1\}^{\omega}$. Then $w\in \EL$ if and only if for all $i\geq 1$, $\pref_w(i) \in \EL$.
\end{enumerate}
\end{lemma}

\begin{definition}[Minimum density, minimum-density prefix, slope] 
Let $w\in \{0,1\}^*\cup \{0,1\}^{\omega}$. 
Define the {\em minimum density of $w$} as $\delta(w) = \inf \{ D_w(i) \mid 1\leq i \}$. If this infimum is attained somewhere, then we also define
$\iota(w) = \min \{ j\geq 1 \mid \forall i: D_w(j) \leq D_w(i) \}$ and $\kappa(w) = P_w(\iota(w)).$
We refer to $\pref_w(\iota(w))$ as the {\em minimum-density prefix}, the shortest prefix with density $\delta(w)$.
For an infinite word $w$, we define the {\em slope} of $w$ as $\lim_{i  \to \infty} D_w(i)$, if this limit exists.
\end{definition}

\begin{remark} Note that $\iota(w)$ is always defined for finite words, while for infinite words, a prefix which attains the infimum may or may not exist.
We note further that density and slope of infinite binary words do not necessarily coincide. In particular, while $\delta(w)$ exists for every $w$, the limit
$\lim_{i  \to \infty} D_w(i)$ may not exist, i.e., $w$ may or may not have a slope.
As an example, consider the word $w = v_0 v_1 v_2 \cdots$, where for each $i$,  $v_i = 1^{2^i} 0^{2^i}.$ Then,
$\delta(w) = 1/2$ and $\lim_{i  \to \infty} D_w(i)$ does not exist, since $D_w(i)$ has an infinite subsequence which is constant $1/2$, and another which tends to $2/3$.

Moreover, even for words $w$ for which the slope is defined, this can be different from the minimum density. If $w$ has slope $\alpha$, then $\alpha = \delta(w)$ if and only if for all $i$, $D_w(i) \geq \alpha$.
For instance, the infinite word $01^\omega$ has slope $1$ but its minimum density is $0$. On the other hand, the infinite word $1(10)^\omega$ has both slope and minimum density $1/2$.
\end{remark}


\section{Operations generating infinite prefix normal words}\label{sec:properties_of_infinite_pnw}

In \cite{CLR18}, we introduced an operation which takes a finite prefix normal word $w$ ending in $1$ and extends it by a run of $0$s followed by a new $1$, in such a way that this new $1$ is placed in the first possible position without violating prefix normality. This operation, called \flipext, leaves the minimum density invariant. Moreover, by repeatedly applying the \flipext\ operation, an infinite prefix normal word is produced which is the densest among all prefix normal words with given prefix $w$. 

Here we extend the definition of \flipext\ to {\em all } prefix normal words containing at least one $1$ and show that the same properties hold, even if the original word $w$ does not end in $1$.

\begin{definition}[Operation $\flipext$]\label{def:flipext}
Let $w\in {\cal L}_{\fin} \setminus \{0\}^*$. Define $\flipext(w)$ as the finite word $w0^k1$, where $k=\min \{ j \mid w0^j1 \in {\cal L}\}$. We further define the infinite word $v = \flipext^{\omega}(w)$. 
\end{definition}

The next proposition is a slightly more general form of Lemma 13 from~\cite{CLR18}: 

\begin{proposition}\label{flipext-properties}
Let $w \in {\cal L}_{\fin} \setminus\{0\}^*$ and $v \in \flipext^*(w) \cup \{\flipext^{\omega}(w)\}$. 
Then $\delta(v) = \delta(w)$, and, as a consequence, $\iota(v) = \iota(w)$
and $\kappa(v)=\kappa(w)$. Moreover,  $D_{v} (j \cdot \iota(w)) = \delta(w)$ for all $j \geq 1$.
\end{proposition}

\begin{proof}
Let $w\in {\cal L}$. If the last character of $w$ is a $1$, then the claim holds by Lemma 13 of~\cite{CLR18}. 

Else $w$ ends in a run of $0$s. Let $\ell$ be the length of this run, and $w'$ be such that $w = w'0^{\ell}$. Let $w'' = \flipext(w') = w'0^k1$, i.e.\ by definition of \flipext, $k$ is minimal s.t.\ $w'0^k1 \in {\cal L}$. If $\ell\leq k$, then $\flipext(w) = \flipext(w') = w''$. Since $w'$ is a prefix of $w$, and $w$ is a prefix of $w''$, we have $\delta(w') \geq \delta(w) \geq \delta(w'')$. Since $w'$ ends in a $1$, $\delta(w'') = \delta(w')$, and thus $\delta(w'') = \delta(w)$. 

Otherwise $\ell>k$, therefore $\flipext(w') = w'0^{\ell'}1 \in {\cal L}_{\fin}$ for some $\ell' < \ell,$ hence $w'0^{\ell}1 \in {\cal L}_{\fin}.$
The latter implies $\flipext(w) = w1$ and $\delta(\flipext(w)) = \delta(w).$

Further iterations $\flipext^{(i)}(w)$ fulfil the claim due to the fact that $\flipext(w)$ ends in a $1$.

We now show the second statement: $D_{v} (j \cdot \iota(w)) = \delta(w)$ for all $j \geq 1$. 
We show it by induction. It is clearly true for 
$j = 1$, moreover for each $j > 1$ assuming $D_{v} ((j-1) \cdot \iota(w)) = \delta(w)$  and letting 
$w' = \pref_v((j-1)\cdot \iota(w))$ and $w''$ be the substring of size $\iota(w)$ such that $w'w'' =  \pref_v(j \cdot \iota(w)),$ we have 

\begin{align*}
\delta(w) &= \delta(v) \leq D_{v} (j \cdot \iota(w)) = \frac{|w'|_1 + |w''|_1}{j \cdot \iota(w)} \\
&\leq 
\frac{P_w(\iota(w)) (j-1) + P_w(\iota(w))}{j \cdot \iota(w)} = \frac{P_w(\iota(w))}{\iota(w)} = \delta(w),
\end{align*}

where in the second inequality  we are using $|w'|_1 = P_w(\iota(w)) (j-1) \iota(w)$ (induction hypothesis) and 
$|w''|_1 \leq P_w(\iota(w))$ (since $v$ is prefix normal).
\qed
\end{proof}

The next proposition states that the infinite word which is generated by repeatedly applying the \flipext\ operation is the densest among all prefix normal words with prefix $w$. 

\begin{proposition}\label{flipext-densest}
	Let $w\in {\cal L}_{\fin}\setminus {0}^*$, $v= \flipext^{\omega}(w)$, and let $z \in {\cal L}_{\inf}$ such that $\pref_{z}(|w|) = w$. Then for every $i = 1, 2, \dots$ we have $P_v(i) \geq P_z(i).$
\end{proposition}
\begin{proof}
We argue by contradiction. Let $i$ be the smallest integer such that $P_v(i) < P_z(i).$ Clearly $i > |w|$ and, by the minimality assumption we must have
$P_v(i-1) = P_z(i-1)$ and $v_i = 0, z_i = 1.$ By definition of $\flipext$ there must exist $j < i$ such that $|v_{j+1} \dots v_{i-1} 1|_1 > P_v(i-j) \geq P_z(i-j),$
for otherwise we would have $v_i = 1.$ Since $v$ is prefix normal, it also follows that we have $|v_{j+1} \dots v_{i-1} v_i|_1 = P_v(i-j) \geq P_z(i-j).$

From this, since by the minimality of $i$ it holds that $P_z(j) \leq P_v(j)$, we have that 
$|z_{j+1} \dots z_{i-1} z_i|_1 = P_z(i) - P_z(j) > P_v(i) - P_v(j) = P_v(i-j) \geq P_z(i-j),$ violating the prefix normality of $z.$ \qed
\end{proof}

We now define a different operation, called \lazyaflipext, which, given a prefix normal word $w$, extends it by adding $0$s as long as the minimum density of the resulting word is not
smaller than $\delta(w)$, and only then adding a $1$. We show that this operation preserves the prefix normality of the resulting word.

\begin{definition}[Operation $\lazyaflipext$]
	Let $\alpha \in (0,1]$ and let $w \in \mathcal{L}_{\fin}$ with $\delta(w)\geq \alpha$. We define $\lazyaflipext(w,\alpha)$ as the finite word $w0^k1$ where $k = \max\{j \mid \delta(w0^j) \geq \alpha\}$. We further define the infinite word $v = \lazyaflipext^{\omega}(w,\alpha)$.
\end{definition}

\begin{example}
	Let $w = 111$ and let $\alpha = \sqrt{2}-1$. Then $\lazyaflipext(w,\alpha) = 11100001$, since $\delta(1110000) = 3/7 \geq \alpha$ and $\delta(11100000) = 3/8 < \alpha$. Further, $\lazyaflipext^{(2)}(w,\alpha) = 1110000101$, since $\delta(111000010) = 4/9 \geq \alpha$ and $\delta(1110000100) = 2/5 < \alpha$.
\end{example}

\begin{lemma}\label{lemma:lazy is pn}
Let $\alpha \in (0,1].$ For every $w \in \L_{\fin}$ with $\delta(w)\geq \alpha$, the word $v = \lazyaflipext(w,\alpha)$ is also prefix normal, with $\delta(v) \geq \alpha.$
\end{lemma}

\begin{proof} First note that $\delta(v) \geq \alpha$ by definition.
Now write $v = w0^k1$, and let $u = \flipext(w) = w0^\ell 1$. Recall that $\ell = \min\{j \mid w0^j1 \in \EL \}$.  If $k<\ell$, this implies $\delta(u)<\alpha$, in contradiction to Proposition~\ref{flipext-properties}, since $\delta(u)=\delta(w)\geq \alpha$. Thus $k\geq \ell$, from which follows $v\in \EL$. \qed

\end{proof}

\begin{corollary}\label{lemma:lazy of w is pn}
Let $\alpha \in (0,1]$ and $w \in \mathcal{L}_{\rm fin}$ with $\delta(w)\geq \alpha$. Then $v = \lazyaflipext^{\omega}(w,\alpha)$ is an infinite prefix normal
word and $\delta(v) = \alpha$.
\end{corollary}

\begin{proof}
	That $v$ is prefix normal follows from Lemma~\ref{lemma:simple pn facts} and from Lemma~\ref{lemma:lazy is pn}, which also implies that $\delta(v)\geq \alpha$.
	However, if $\delta(v) > \alpha$ was true, then for a suitably long prefix $i$, we would get a contradition to the definition of the $\lazyaflipext$ operation. 	\qed
\end{proof}

Fix $w \in {\cal L}_{\fin}$. The next proposition states that the \lazyflipext\ operation with  $\alpha = \delta(w)$, applied to $w$,  generates a prefix normal word that has the minimum number of $1$s among all prefix normal words with prefix $w$ 
and  minimum density $\delta(w)$. 

\begin{proposition}\label{lazy-alpha-flipext-least-dense}
	Let $w\in {\cal L}_{\fin}$, $\alpha = \delta(w)$, $v= \lazyaflipext^{\omega}(w,\alpha)$, and $z \in {\cal L}_{\inf}$ such that $\pref_{z}(|w|) = w$ and $\delta(z) \geq \delta(w)$. Then for all $i = 1, 2, \dots$, we have $P_v(i) \leq P_z(i).$ 
	\end{proposition}

\begin{proof} We argue by contradiction. Let $i$ be the smallest integer such that $P_v(i) < P_z(i).$ Clearly $i > |w|$ and, by the minimality assumption, we have
$P_v(i-1) = P_z(i-1)$ and $v_i = 0, z_i = 1.$ Let $u = \pref_v(i-1)$. Since $i>|w|$ and $v_i=1$, therefore $u1 = \lazyaflipext(u', \alpha)$ for some $u'$, and thus, by definition of $\lazyaflipext$, $P_{u0}(i)/i < \alpha$. But $u0 = \pref_i(z)$, so we have 

$$\delta(z) \leq D_z(i) = \frac{P_z(i)}{i} = \frac{P_{u0}(i)}{i} < \delta(w),$$

in contradiction to the density of $z$. \qed

\end{proof}

\begin{theorem}\label{thm:extremal}
Let $w\in \EL_{\fin}$ with $\alpha = \delta(w) \in (0,1]$, and let $z \in \EL_{\inf}$ such that $\pref_{z}(|w|) = w$ and $\delta(z)\geq \alpha$. Let $u = \flipext^{\omega}(w)$ and $v = \lazyaflipext^{\omega}(w,\alpha)$. Then $v \leq_{\lex} z \leq_{\lex} u$. 
\end{theorem}

\begin{proof}
Follows from Prop.~\ref{flipext-densest} and  Prop.~\ref{lazy-alpha-flipext-least-dense}. \qed
\end{proof}

Note that if $\pref_z(|w|)=w$, then $\delta(z) \geq \delta(w)$ implies that, in fact, $\delta(z) = \delta(w)$ holds, since $z$ is an extension of $w$. Theorem~\ref{thm:extremal} states then that all prefix normal extensions of $w$ with the same minimum density as $w$ lie lexicographically between the  \lazyaflipext- and the \flipext-extensions of $w$. However, not all extensions of $w$ between these two words are prefix normal, as we can see in the next example.

\begin{example}\label{ex:flipext_vs_lazyalphaflipext}
Let $w = 1101101100100010000001$, with $\alpha = \delta(w) = 8/21$, then 
\begin{eqnarray*}
	v &=& \textstyle{\lazyaflipext^{(8)}}(w,\alpha) = w01001010010010100100,\\
	u &=& \textstyle{\flipext^{(8)}(w)} = w101101100100010000001.
\end{eqnarray*}
	
	\medskip
	
Let $p = w100111010100000100001$ and $q = w101101010100001000001$, we have that for all $1 \leq i \leq 42$, $P_v(i) \leq P_p(i), P_q(i) \leq P_u(i)$ and $v \leq_\lex p, q \leq_\lex u$. Note that $p$ is not prefix normal, while $q$ is prefix normal. 

\begin{figure}[!htb]
\begin{center}
\includegraphics[width=.95\textwidth]{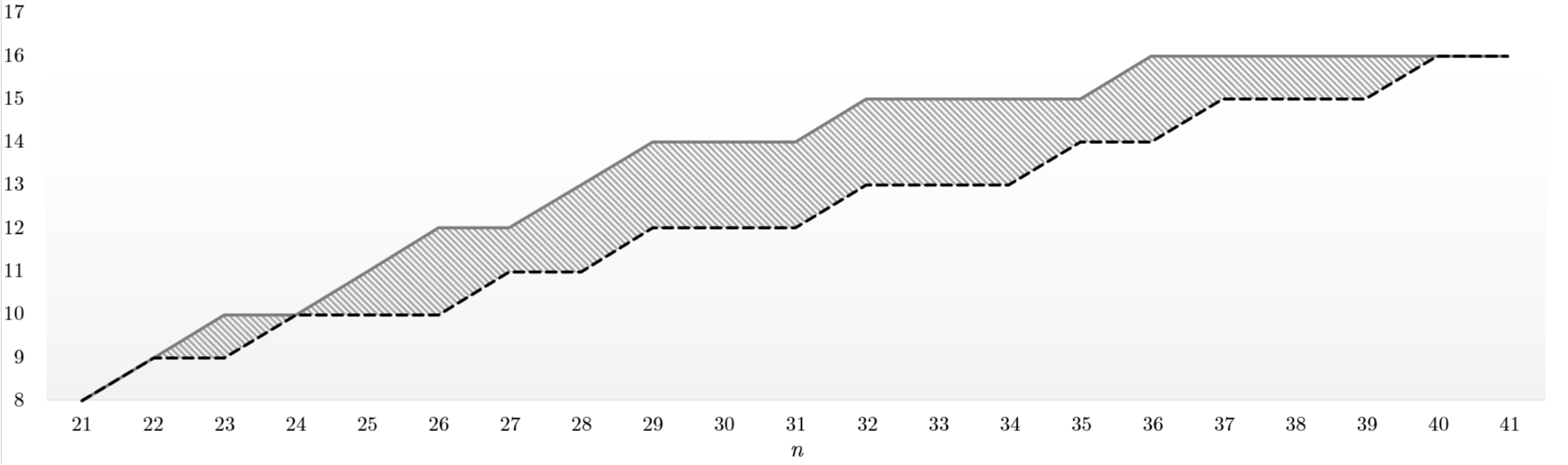}
\caption{Given $w = 1101101100100010000001$ the plot represents the last characters of $\flipext^{(8)}(w)$ (solid) and the $\lazyaflipext^{(8)}(w,\alpha)$ (dashed). See Example~\ref{ex:flipext_vs_lazyalphaflipext}. A $1$ corresponds to a diagonal segment in direction NE, while a $0$ to one in direction SE. On the $x$-axis we have the length of the prefix, and on the $y$-axis, the number of $1$s minus the number of $0$s  in the prefix.
The shaded area contains all prefix normal words with $w$ as prefix and minimum density equal to $\delta(w)$. Note, however, that not all words in that area are prefix normal. \label{fig:flipext_vs_lazy}}
\end{center}
\vspace{-0.7cm}
\end{figure}
\vspace{-.2cm}

\end{example}


\section{Sturmian words and prefix normal words}\label{sec:sturmian}

In the previous section, we presented operations that construct infinite prefix normal words by extending finite prefix normal words. In particular, the $\lazyaflipext$ operation extends a finite binary word with as few $1$s as possible while preserving its minimum density. This is reminiscent of the characterization of Sturmian words in terms of mechanical words and the slope.
Led by this analogy, in this section we provide a complete characterization of  Sturmian words which are prefix normal.
We refer the interested reader to~\cite[Chapter 2]{Lothaire3}, for a comprehensive treatment of
Sturmian words. Here we briefly recall some facts which we will need later.

\begin{definition}[Sturmian words]
Let $w\in \{0,1\}^{\omega}$. Then $w$ is called {\em Sturmian} if it is balanced and aperiodic.
\end{definition}

An equivalent definition of Sturmian words is that they are irrational mechanical, a definition we recall next.

\begin{definition}[Mechanical words]
Given two real numbers $0\leq \alpha \leq 1$ and $0 \leq \tau < 1$, the {\em lower mechanical word} $s_{\alpha,\tau} = s_{\alpha,\tau}(1)\, s_{\alpha,\tau}(2)\, \cdots $
and  the {\em upper mechanical word} $s'_{\alpha,\tau} = s'_{\alpha,\tau}(1)\, s'_{\alpha,\tau}(2)\, \cdots$ are given by

\[\begin{aligned}
s_{\alpha,\tau} (n) = \lfloor \alpha n + \tau\rfloor - \lfloor \alpha (n-1)  + \tau\rfloor \\
s'_{\alpha,\tau} (n) = \lceil \alpha n  + \tau\rceil - \lceil \alpha (n-1) + \tau\rceil
\end{aligned}\quad\quad (n\geq 1).\]

Then $\alpha$ is called the {\em slope} and $\tau$ the {\em intercept} of $s_{\alpha, \tau}, s'_{\alpha, \tau}$.
A word $w$ is called {\em mechanical} if $w = s_{\alpha,\tau}$
or $w=s'_{\alpha,\tau}$ for some $\alpha,\tau$.
It is called {\em rational mechanical} (resp.\ {\em irrational mechanical}) if $\alpha$ is rational (resp.\ irrational).
\end{definition}

\begin{fact}[Some facts about Sturmian words~\cite{Lothaire3}]\label{fact:sturmian-properties}
\begin{enumerate}
\item An infinite binary word is Sturmian if and only if it is irrational mechanical.
\item For $\tau=0$ and irrational $\alpha$, there exists a word $c_{\alpha}$, called the {\em characteristic word with slope $\alpha$}, s.t.\  $s_{\alpha,0} = 0c_{\alpha}$ and $s'_{\alpha,0} = 1c_{\alpha}$. This word $c_{\alpha}$ is a Sturmian word itself, with both slope and intercept $\alpha$.
\item For two Sturmian words $w$ and $v$ with the same slope, $\Fac(w) = \Fac(v)$.
\end{enumerate}
\end{fact}

We now show that the word  $\lazyaflipext^{\omega}(1,\alpha)$ coincides with the  upper mechanical word $s'_{\alpha, 0}$. This also implies that $s'_{\alpha, 0}$ is prefix normal, as noted in the subsequent corollary.

\begin{lemma}\label{lemma:mechanical equals lazy}
Fix $\alpha \in  (0,1]$ and let $v = \lazyaflipext^{\omega}(1,\alpha)$. Let $s = s'_{\alpha,0}$ be the upper mechanical word of slope $\alpha$ and intercept
$0$. Then $v = s$.
\end{lemma}

\begin{proof}
Let $s_i$ and $v_i$ denote the $i$th character of $s$ and $v$ respectively. We argue by induction on $i$ that $v_i = s_i$. The claim is true for
$i=1$ since, directly from the definitions we have $v_1 = 1 = s_1.$
Let $n > 1$ and assume that for each $i < n$ we have $v_i = s_i.$ For the induction step we argue according to the character $s_n.$

(i) If $s_n = 1$, by definition $\lceil n \alpha \rceil - \lceil (n-1) \alpha\rceil = 1.$ Thus, $\lceil (n-1)\alpha \rceil < n\alpha.$
Using this inequality and the induction hypothesis together with the
definition of $s'_{\alpha, 0}$ we have that  $|v_1 \cdots v_{n-1}|_1 = |s_1 \cdots s_{n-1}|_1 = \lceil (n-1) \alpha \rceil < \alpha n.$
Therefore $|v_1 \cdots v_{n-1} 0|_1 = |v_1 \cdots v_{n-1}|_1  < \alpha n$, which means that $\delta(v_1 \cdots v_{n-1} 0)$ $< \alpha,$ hence
by definition $\lazyaflipext(v_1\cdots v_{n-1},\alpha) = v_1 \cdots v_{n-1} 1$, i.e., $v_n = 1 = s_n.$

(ii) If $s_n = 0$, by definition $\lceil n \alpha \rceil - \lceil (n-1) \alpha\rceil = 0.$ Thus, $\lceil (n-1)\alpha \rceil \geq n\alpha.$
Using this inequality and the induction hypothesis together with the
definition of $s'_{\alpha, 0}$ we have that  $|v_1 \cdots v_{n-1}|_1 = |s_1 \cdots s_{n-1}|_1 = \lceil (n-1) \alpha \rceil \geq \alpha n.$
Therefore $|v_1 \cdots v_{n-1} 0|_1 = |v_1 \cdots v_{n-1}|_1  \geq \alpha n$ which means that $\delta(v_1 \cdots v_{n-1} 0)$ $\geq \alpha,$ hence
by definition $\lazyaflipext(v_1\cdots v_{n-1},\alpha) = v_1 \cdots v_{n-1} 0\cdots01$, i.e., $v_n = 0 = s_n.$
	\qed
\end{proof}

\begin{corollary}\label{lemma:mechanical is pn}
	Let $\alpha \in (0,1]$. Then $s'_{\alpha,0}$ is an infinite prefix normal word and $\delta(s'_{\alpha,0}) = \alpha$.
\end{corollary}
The following theorem fully characterizes those Sturmian words which are prefix normal.
\begin{theorem}\label{thm:characteristic}
A Sturmian word $s$ of slope $\alpha$
is prefix normal if and only if $s = 1c_{\alpha}$, where $c_{\alpha}$ is the characteristic Sturmian
word with slope $\alpha$.
\end{theorem}
\begin{proof}
By definition, $\alpha$ is irrational. Let $s = s'_{\alpha,0}.$ Then $s$ is Sturmian and prefix normal by Corollary  \ref{lemma:mechanical is pn}.
Let $t$ be a Sturmian word with the same slope $\alpha$ which is also prefix normal. By Fact~\ref{fact:sturmian-properties}, $s$ and $t$ have the same factors.

Assume, by contradiction, that $s \neq t$, hence there exists $i \geq 1$ such that $|s_1\cdots s_i|_1 \neq |t_1\cdots t_i|_1.$
Assume,  without loss of generality (since we can, if necessary, swap $s$ and $t$ in the following argument),
that $|s_1 \cdots s_i|_1 > |t_1 \cdots t_i|_1.$ Then, since $s_1 \cdots s_i$ is also a factor of $t$, there is a $j \geq1$ such that $t_{j+1} \cdots t_{j+i} =
s_1 \cdots s_i$, hence   $|t_{j+1} \cdots t_{j+i}|_1 > |t_1 \cdots t_i|_1$ contradicting the assumption that $t$ is prefix normal.
	\qed
\end{proof}



\section{Prefix normal words, prefix normal forms,  and abelian complexity}\label{sec:pnf}

Given an infinite word $w$, the {\em abelian complexity} function of $w$, denoted $\psi_w$, is given by $\psi_w(n) = |\{ \pv(u) \mid u\in\Fac(w), |u|=n\}|$, the number of Parikh vectors of $n$-length factors of $w$. A word $w$ is said to have bounded abelian complexity if there exists a $c$ s.t.\ for all $n$, $\psi_w(n) \leq c$. Note that a binary word is $c$-balanced if and only if its abelian complexity is bounded by $c+1$. We denote the set of Parikh vectors of factors of a word $w$ by $\Pi(w) = \{ \pv(u) \mid u \in \Fac(w)\}$. Thus, $\psi_w(n) = |\Pi(w) \cap \{(x,y) \mid x+y = n\}|$.
In this section, we study the connection between prefix normal words and abelian complexity.

\subsection{Balanced and $c$-balanced words.}

Based on the examples in the introduction, one could conclude that any word with bounded abelian complexity can be turned into a prefix normal word by prepending a fixed number of $1$s. However, consider the word $w=01^{\omega}$, which is balanced, i.e.\ its abelian complexity function is bounded by $2$. It is easy to see that $1^k w \not\in \EL$ for every $k\in \IN$.

Sturmian words are precisely the words which are aperiodic and whose abelian complexity is constant $2$~\cite{RSZ11}. For Sturmian words, it is always possible to prepend a finite number of $1$s to get a prefix normal word, as we will see next. Recall that for a Sturmian word $w$, at least one of $0w$ and $1w$ is Sturmian, with both being Sturmian if and only if $w$ is characteristic~\cite{Lothaire3}.

\begin{lemma}
Let $w$ be a Sturmian word with slope $\alpha$. Then

\begin{enumerate}
\item
$1w \in \EL$ if and only if $0w$ is Sturmian,
\item
if $0w$ is not Sturmian, then $1^{n}w \in \EL$ for $n = \lceil 1/(1-\alpha) \rceil$.
\end{enumerate}
\end{lemma}

\begin{proof}
{\em 1.} Let $0w$ be Sturmian and let $u$ be some factor of $1w$. If $u$ is a prefix of $1w$, there is nothing to show, therefore let $u\in\Fac(w)$, with $|u|=n$ and $|u|_1=k$. Since $0w$ is Sturmian, we have that the prefix of $0w$ of length $n$ has at least $k-1$ $1$s, thus $P_{1w}(n) \geq k = |u|_1$, as desired. Conversely, if $0w$ is not Sturmian, this means that it is not balanced, therefore there exists a factor $u$ of $w$ s.t.\ $||u|_1 - |0w_1\cdots w_{n-1}|_1| \geq 2$, where $|u|=n$. Since $w$ is Sturmian, we have that $||w_1\cdots w_{n-1}|_1 - |u_1\cdots u_{n-1}|_1| \leq 1$ and $||w_1\cdots w_{n-1}|_1 - |u_2\cdots u_{n}|_1| \leq 1$. Let $|w_1\cdots w_{n-1}|_1 = k$, then this implies, by a case-by-case consideration, that $|u_1\cdots u_{n-1}|_1 = |u_2\cdots u_{n}|_1= k+1$, and thus $|1w_1\cdots w_{n-1}|_1 = k +1 < k+2 = |u|_1$, showing that $1w$ is not prefix normal.

{\em 2.} First note that a Sturmian word of slope $\alpha$ cannot have a run of $1$s of length 
$\lceil 1/(1-\alpha)\rceil.$  
To see this, it is enough to
consider the upper mechanical word of slope $\alpha$ and intercept $0$ (since all the other words with the same slope have the same set of factors).
Let us write $s = s_{\alpha, 0} = s_1 s_2 \cdots$

Now $s$ has a run of $n$ 1s if and only if there exists an $i \geq 0$ such that $s_{i+1} = s_{i+2} = \cdots = s_{i+n} = 1.$
By the definition of mechanical words, we have that the last condition is equivalent to
$$\lceil \alpha(i+n) \rceil - \lceil \alpha i \rceil = n.$$

On the other hand, if $n \geq \frac{1}{1-\alpha}$, i.e., $\alpha \leq  \frac{n-1}{n}$ we have that the sum of the character
$\sum_{j=1}^n s_{i+j}$ satisfies

\begin{eqnarray*}
\sum_{j=1}^n s_{i+j}  &= & \lceil \alpha(i+n) \rceil - \lceil \alpha i \rceil 
\leq \lceil \alpha i \rceil + \lceil \alpha n \rceil - \lceil \alpha i \rceil \\
&=& \lceil \alpha n \rceil 
<  \alpha n + 1 \leq  \frac{n-1}{n} \times n + 1 = n.
\end{eqnarray*}

\noindent i.e., strictly smaller than $n$, i.e., we have a contradiction $s_{i+1} \cdots s_{i+n} \neq 1^n$.

Now fix $n = \lceil 1/(1-\alpha) \rceil$ and  let  $w' = 1^{n}w$.  Let  $u\in\Fac(w)$. Since, as shown above, $1^n$ is not a factor, if $|u|\leq n$, there is nothing to show. So let $|u| = n+m$. Then $|u_1\cdots u_n|_1 \leq n-1$, and since $w$ is balanced, we have that $|w_1\cdots w_{m}|_1 \geq |u_{n+1}\cdots u_{n+m}|_1 - 1$, yielding that $P_{w'}(n+m) \geq n + |u_{n+1}\cdots u_{n+m}|_1 - 1 \geq |u|_1$. \qed
\end{proof}

\begin{lemma}\label{lemma:c-bal}
Let $w$ be a $c$-balanced word. If there exists a positive integer $n$ s.t.\ $1^n \not\in\Fac(w)$, then the word $z = 1^{nc}w$ is prefix normal.
\end{lemma}
\begin{proof}
We are going to show that every factor  $u$ of $z$ satisfies the prefix normal condition $|u|_1 \leq P_z(|u|)$.
It is not hard to see that we can limit ourselves to only
considering  factors $u$ such that $u$ does not overlap with the prefix of $z$ of the same length.

If $|u| \leq nc$ then $|u|_1 \leq |u| = P_z(|u|).$ Assume now that $u = u' u''$ with $|u'| = nc$ and $|u''| > 0.$
Since $u'$ is a factor of $w$ of size $nc$ the condition that $w$ does not contain a factor $1^n$ implies that
$u'$ contains at least $c$ 0s, i.e., $|u'|_1 \leq |u'| - c.$
Moreover, since $w$ is $c$-balanced, we have that $|u''|_1 \leq P_w(|u''|)+c.$
Therefore, observing that $\pref_z(|u|) = \pref_z(|u'| + |u''|) = 1^{nc} \pref_w(|u''|)$ we have that
$P_z(|u|) = nc + P_w(|u''|) \geq |u'|_1 + |u''|_1 = |u|_1.$ 
\qed
\end{proof}

In particular, Lemma~\ref{lemma:c-bal} implies that any $c$-balanced word with infinitely many $0$s can be turned into a prefix normal word by prepending a finite number of $1$s, since such a word cannot have arbitrarily long runs of $1$s.
Note, however, that the number of $1$s to prepend from Lemma~\ref{lemma:c-bal}  is not tight, as can be seen e.g.\ from the Thue-Morse word $\textswab{t}$: the longest run of $1$s in \textswab{t} is $2$ and \textswab{t} is $2$-balanced, but $11\textswab{t}$ is prefix normal, as will be shown in the next section (Lemma~\ref{lemma:11tm}).

\subsection{Prefix normal forms and abelian complexity.}

Recall that for a word $w$, $F^a_w(i)$ is the maximum number of $a$'s in a factor of $w$ of length $i$, for $a\in\{0,1\}$.

\begin{definition}[Prefix normal forms]
Let $w\in \{0,1\}^{\omega}$. Define the words $w'$ and $w''$ by setting, for $n\geq 1$, $w'_n = F^1_w(n) - F^1_w(n-1)$ and $w''_n = \overline{F^0_w(n) - F^0_w(n-1)}$.
We refer to $w'$ as the {\em prefix normal form of $w$ w.r.t.\ $1$} and to $w''$ as the {\em prefix normal form of $w$ w.r.t.\ $0$}, denoted $\PNF_1(w)$ resp.\ $\PNF_0(w)$.
\end{definition}

In other words, $\PNF_1(w)$ is the sequence of first differences of the maximum-$1$s function $F^1_w$ of $w$. Similarly, $\PNF_0(w)$ can be obtained by complementing the sequence of first differences of the maximum-$0$s function $F^0_w$ of $w$. Note that for all $n$ and $a\in\{0,1\}$, either $F^a_w(n+1) = F^a_w(n)$ or $F^a_w(n+1) = F^a_w(n) + 1$,  and therefore $w'$ and $w''$ are words over the alphabet $\{0,1\}$. In particular, by construction, the two prefix normal words allow us to recover the maximum-$1$s and minimum-$1$s functions of $w$:

\begin{observation}\label{obs1} Let $w$ be an infinite binary word and $w' = \PNF_1(w), w''=\PNF_0(w)$. Then
$P_{w'}(n) =  F^1_w(n)$ and $P_{w''}(n) =  n-F^0_w(n) = f^1_w(n).$
\end{observation}

\begin{lemma}
Let $w\in \{0,1\}^{\omega}$. Then $\PNF_1(w)$ is the unique $1$-prefix normal word $w'$ s.t.\ for all $i\in \IN$, $F^1_{w'}(i) = F^1_w(i)$. Similarly, $\PNF_0(w)$ is the unique $0$-prefix normal word $w''$ s.t.\ for all $i\in \IN$, $F^0_{w''}(i) = F^0_w(i)$.
\end{lemma}

\begin{proof}
Let $w' = \PNF_1(w)$ and $w'' = \PNF_0(w)$. First note that, by construction, for all $i\in \IN$, $F^1_{w'}(i) = F^1_w(i)$ and $F^0_{w''}(i) = F^0_w(i)$. It is easy to see that $w'$ is $1$-prefix normal and $w''$ is $0$-prefix normal. For uniqueness, note that for $a\in \{0,1\}$ and an $a$-prefix normal word $v$, we have $\PNF_a(v) = v$. \qed
\end{proof}

\begin{example}
The two prefix normal forms and the maximum-$1$s and maximum-$0$s functions of the Fibonacci word $\textswab{f} = 01001010010010100101\cdots$ are given in Table~\ref{tab:fib}.

\begin{table}[h]
{\small
    \begin{tabularx}{\textwidth}{c|YYYYYYYYYYYYYYYYYYYY}
        \toprule
        n & 1 & 2 & 3 & 4 & 5 & 6 & 7 & 8 & 9 & 10 & 11 & 12 & 13 & 14 & 15 & 16 & 17 & 18 & 19 & 20\\
        \cmidrule(l){1-21}%
        $F^0_\textswab{f}(n)$ & 1 & 2 & 2 & 3 & 4 & 4 & 5 & 5 & 6 & 7 & 7 & 8 & 9 & 9 & 10 & 10 & 11 & 12 & 12 & 13 \\
        $F^1_\textswab{f}(n)$ & 1 & 1 & 2 & 2 & 2 & 3 & 3 & 4 & 4 & 4 & 5 & 5 & 5 & 6 & 6 & 7 & 7 & 7 & 8 & 8 \\
        \cmidrule(l){1-21}%
        $\PNF_0(\textswab{f})$ & 0 & 0 & 1 & 0 & 0 & 1 & 0 & 1 & 0 & 0 & 1 & 0 & 0 & 1 & 0 & 1 & 0 & 0 & 1 & 0 \\
        $\PNF_1(\textswab{f})$ & 1 & 0 & 1 & 0 & 0 & 1 & 0 & 1 & 0 & 0 & 1 & 0 & 0 & 1 & 0 & 1 & 0 & 0 & 1 & 0 \\
        \bottomrule
    \end{tabularx}%
    \caption{The maximum number of $0$s and $1$s ($F^0_\textswab{f}(n)$ and $F^1_\textswab{f}(n)$ resp.) for all $n = 1, \ldots, 20$ of the Fibonacci word $\textswab{f}$, and the prefix normal forms of \textswab{f}.}
    \label{tab:fib}
    }
\end{table}

\end{example}
\begin{figure}
\begin{center}
\includegraphics[width=.95\textwidth]{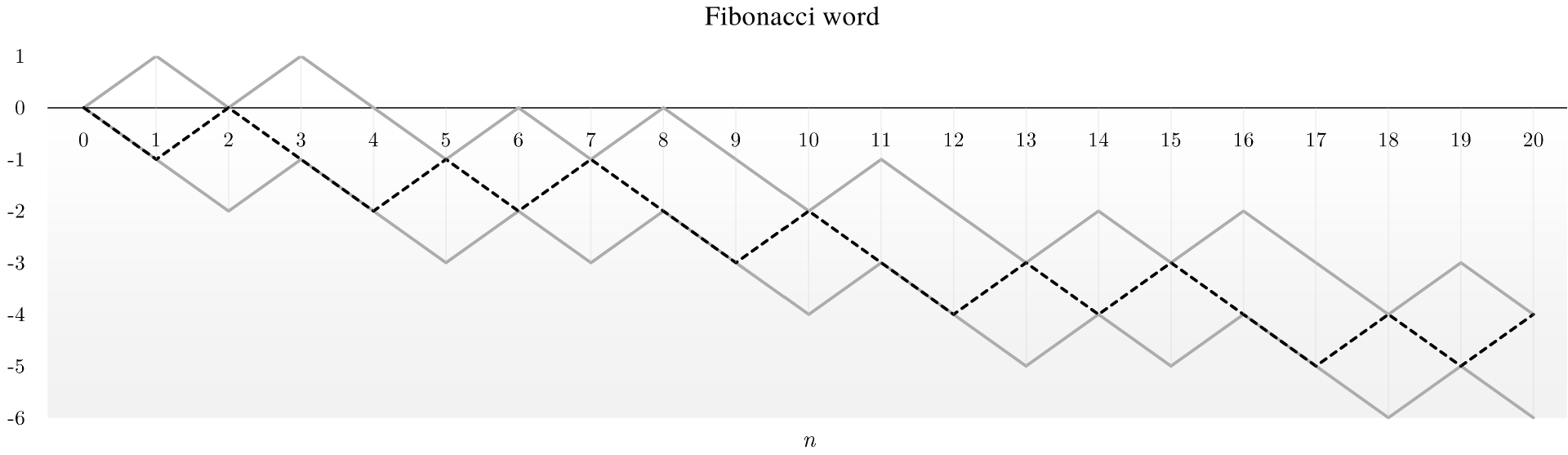}
\caption{The Fibonacci word (dashed) and its prefix normal forms (solid). 
  \label{fig:Fibonacci}}
\end{center}
\vspace{-0.7cm}
\end{figure}
\vspace{-.2cm}

Now we can connect the prefix normal forms of $w$ to the abelian complexity of $w$ in the following way. Given $w' = \PNF_1(w)$ and $w''=\PNF_0(w)$, the number of Parikh vectors of $k$-length factors is precisely $1$ more than the difference in $1$s in the prefix of length $k$ of $w'$ and of $w''$. For example, Fig.~\ref{fig:Fibonacci} shows the prefix normal forms of the Fibonacci word.
The vertical line at 5 cuts through points $(5,-1)$ and $(5,-3)$: the first component stands for the length of the string, the second for the difference between the number of 0s and the number of 1s, therefore indicating Parikh vectors $(2,3)$ and $(1,4)$. 

\begin{figure}
\begin{center}
\includegraphics[width=0.95\textwidth]{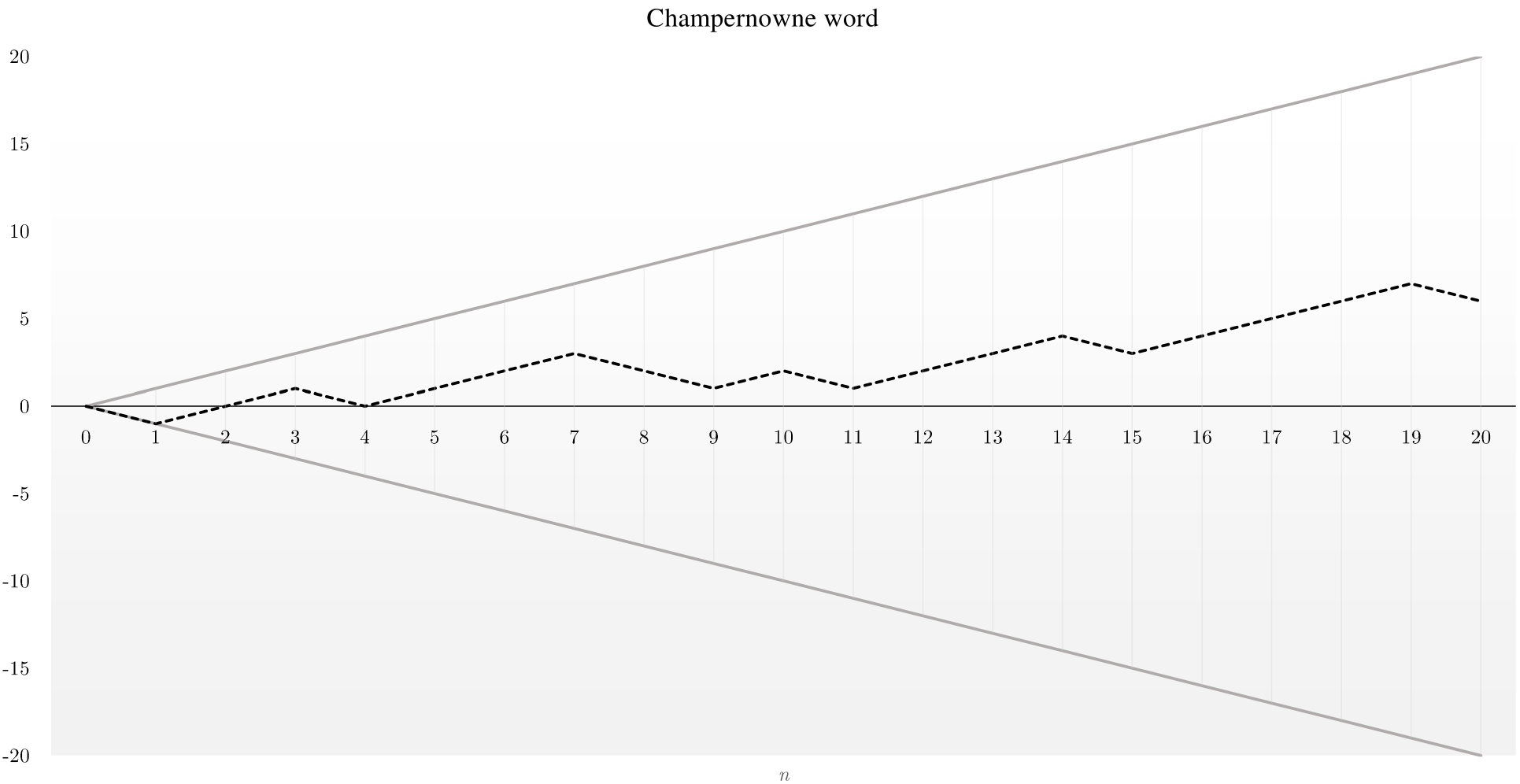}
\caption{The Champernowne word (dashed) and its prefix normal forms (solid). \label{fig:Champernowne}}
\end{center}
\end{figure}

The Fibonacci word, being a Sturmian word, has constant abelian complexity $2$. An example of  a word with unbounded abelian complexity is the Champernowne word, whose prefix normal forms are $1^{\omega}$ resp.\ $0^{\omega}$. (Fig.~\ref{fig:Champernowne}).

\begin{theorem}\label{thm:pnfs_abcomp}
Let $w,v\in \{0,1\}^\omega$. 
\begin{enumerate}
\item $\psi_w(n) = P_{w'}(n) - P_{w''}(n) +1,$ where $w' = \PNF_1(w)$ and $w'' = \PNF_0(w)$.
\item $\Pi(w) = \Pi(v)$ if and only if $\PNF_0(w) = \PNF_0(v)$ and $\PNF_1(w) = \PNF_1(v)$. %
\end{enumerate}
\end{theorem}

\begin{proof}

{\em 1.} Fix an integer $n\geq 1$. By definition, we have that for every factor $u$ of $w$ of length $n$ we have
$n-F^0_w(n) \leq |u|_1 \leq F^1_w(n).$ Therefore $\psi_w(n) \leq F^1_w(n) - (n - F^0_w(n))+1.$

Conversely, since $w$ contains a factor $u'$ of length $n$ with $F^1_w(n)$ many 1s and a factor $u''$ of length $n$ with
$n-F^0_w(n)$ many 1s, if we scan $w$ between an occurrence of $u'$ and an occurrence of $u''$, for each $x \in \{|u''|_1, \dots, |u'|_1\}$
there must be a factor $u'''$ of size $n$ such that $|u'''|_1 = x.$ Therefore $\psi_w(n) \geq F^1_w(n) - (n - F^0_w(n))+1.$
We can conclude that $\psi_w(n) = F^1_w(n) - (n - F^0_w(n)) + 1.$ The desired result then follows by observing that
$n-F^0_w(n) = n - |\pref_{\PNF_0(w)}(n)|_0 = P_{\PNF_0(w)}(n)$ and  $F^1_w(n) = P_{\PNF_1(w)}(n).$

\smallskip

{\em 2.} Follows directly from Observation~\ref{obs1}.
\qed
\end{proof}

Theorem~\ref{thm:pnfs_abcomp} implies that if we know the prefix normal forms of a word, then we can compute its abelian complexity. Conversely, the abelian complexity is the {\em width} of the area enclosed by the two words $\PNF_1(w)$ and $\PNF_0(w)$. In general, this fact alone does not give us the PNFs; but if we know more about the word itself, then we may be able to compute the prefix normal forms, as we will see in the case of the paperfolding word.

We will now give two examples of the close connection between abelian complexity and prefix normal forms, using some recent results about the abelian complexity of infinite words.


\subsubsection{The paperfolding word}
The first few characters of the ordinary paperfolding word
are given by 
\[ \textswab{p} = 0010011000110110001001110011011 \cdots \]

The paperfolding word was originally introduced in~\cite{DavisKnuth70}. One definition is given by: $\textswab{p}_n = 0$ if $n' \equiv 1 \bmod 4$ and $\textswab{p}_n = 1$ if $n' \equiv 3 \bmod 4$, where $n'$ is the unique odd integer such that $n = n'2^k$ for some $k$~\cite{MadillR13}. The abelian complexity function of the paperfolding word was fully determined in~\cite{MadillR13}, giving the following initial values for $\psi_{\textswab{p}}(n)$, for $n\geq 1$: $2,3,4,3,4,5,4,3,4,5,6,5,4,5,4,3,4,5,6,5$, and a recursive formula for the computation of all values. The authors note that for the paperfolding word, it holds that if $u\in\Fac(\textswab{p})$, then also $\overline{u^{\rm rev}} \in \Fac(\textswab{p})$. This implies

\[ F^1_{\textswab{p}}(n) = F^0_{\textswab{p}}(n) \text{ for all $n$, and thus } \PNF_0(\textswab{p}) = \overline{\PNF_1(\textswab{p})}.\]

\begin{figure}
\begin{center}
\includegraphics[width=0.95\textwidth]{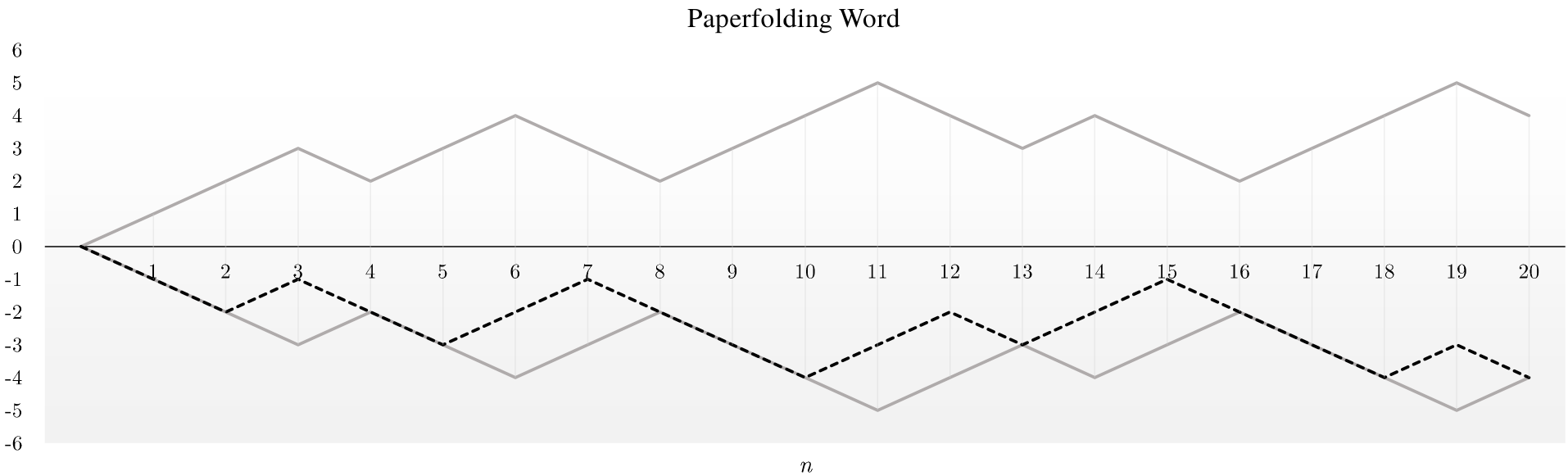}
\caption{The paperfolding word (dashed) and its prefix normal forms (solid). \label{fig:paperfolding}}
\end{center}
\end{figure}

Moreover, from Thm.~\ref{thm:pnfs_abcomp} we get that $F^1_{\textswab{p}}(n) = P_{\PNF_1(\textswab{p})}(n)= (\psi_{\textswab{p}}(n) + n -1)/2$, and thus we can determine the prefix normal forms of {\textswab{p}}, see Fig.~\ref{fig:paperfolding}. 

This same argument holds in general as long as the word has the symmetric property similar to the paperfolding word.
Therefore, we have proved the following lemma.

\begin{lemma} \label{lemma:psi_pnf}
	Let $w\in\{0,1\}^{\omega}$. If for all $u\in \Fac(w)$, it holds that $\overline{u}\in\Fac(w)$ or $\overline{u^{\rm rev}} \in \Fac(w)$, then
$F^1_{w}(n) = F^0_{w}(n) \text{ for all $n$}, \PNF_0(w) = \overline{\PNF_1(w)},$ and $F^1_{w}(n) = (\psi_{w}(n) + n -1)/2.$
\end{lemma}


\subsubsection{Morphic images under the Thue-Morse morphism}
The Thue-Morse word beginning with $0$, which we denote by {\textswab{t}}, is one of the two fix points of the Thue-Morse morphism $\mu_{\tm}$, where $\mu_{\tm}(0) = 01$ and $\mu_{\tm}(1) = 10$:

\[\textswab{t}=\mu_{\tm}^{\omega}(0) = 01101001100101101001011001101001\cdots\]

\begin{figure}
\begin{center}
\includegraphics[width=.95\textwidth]{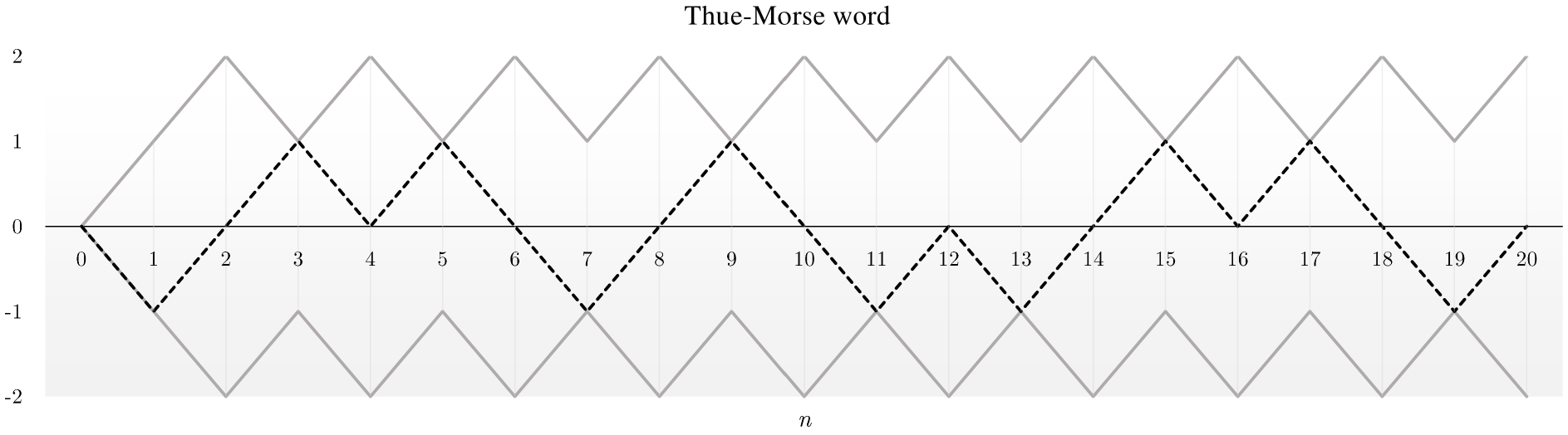}
\caption{The Thue-Morse word (dashed) and its prefix normal forms (solid). \label{fig:Thue-Morse}}
\end{center}
\vspace{-0.5cm}
\end{figure}

The word ${\textswab{t}}$ has abelian complexity function $\psi_{\textswab{t}}(n) = 2$ for $n$ odd and $\psi_{\textswab{t}}(n) = 3$ for $n>1$ even~\cite{RSZ11}. Since {\textswab{t}} fulfils the condition that $u\in\Fac(\textswab{t})$ implies $\overline{u}\in\Fac(\textswab{t})$, we can apply Lemma~\ref{lemma:psi_pnf}, and compute the prefix normal forms of \textswab{t} as $\PNF_1(\textswab{t}) = 1(10)^{\omega}$ and $\PNF_0(\textswab{t}) = 0(01)^{\omega}$, see Fig.~\ref{fig:Thue-Morse}.

For the proof of the abelian complexity of {\textswab{t}} in~\cite{RSZ11}, the Parikh vectors were computed for each length, so we do not really need Lemma~\ref{lemma:psi_pnf} but could have obtained the prefix normal forms directly. Moreover, a much more general result was given in~\cite{RSZ11}:

\begin{theorem}[\hspace{1sp}\cite{RSZ11}]\label{thm:paperfolding}
Let $w$ be an aperiodic infinite binary word. Then $\psi_w = \psi_{\textswab{t}}$ if and only if $w = \mu_{\tm}(w')$, $w = 0\mu_{\tm}(w')$, or $w=1\mu_{\tm}(w')$, for some word $w'$.
\end{theorem}

The abelian complexity function does not in general determine the prefix normal forms, as can be seen on the example of Sturmian words, which all have the same abelian complexity function but different prefix normal forms. However, $\psi_{\textswab{t}}$ does, due to its values $\psi_{\textswab{t}}(n) = 2$ for $n$ odd and $\psi_{\textswab{t}}(n)=3$ for $n$ even, and to the fact that both $F^1_{\textswab{t}}$ and $F^0_{\textswab{t}}$ have difference function with values from $\{0,1\}$:  notice that the only pair of such functions with width $2$ resp.\ $3$ are the PNFs of $\textswab{t}$. Therefore, we can deduce the following from Theorem~\ref{thm:paperfolding}:

\begin{corollary}\label{coro:tm}
For an aperiodic infinite binary word $w$, $\PNF_1(w) = 1(10)^{\omega}$ and $\PNF_0=0(01)^{\omega}$ if and only if $w = \mu_{\tm}(w')$, $w = 0\mu_{\tm}(w')$, or $w=1\mu_{\tm}(w')$, for some word $w'$.
\end{corollary}

To conclude this section, we return to the question of how many $1$s need to be prepended to make the Thue-Morse word prefix normal.

\begin{lemma}\label{lemma:11tm}
	We have $11\textswab{t} \in \EL$. Moreover, this is minimal since $1\textswab{t}$ is not prefix normal.
\end{lemma}

\begin{proof}
	We will show that for every prefix, the number of $1$s in the prefix of $11\textswab{t}$ is greater than or equal to the the number of $1$s in the prefix of $\PNF_1(\textswab{t})$ of the same length. 	Let $v = \PNF_1(\textswab{t})$ and $u=11\textswab{t}$. It is easy to see that $P_v(n) =\lfloor \frac{n}{2} \rfloor + 1$ and $$P_u(n) = \begin{cases}
	\frac{n}{2}+1 & \text{if } n \text{ is even}\\
	\lfloor \frac{n}{2} \rfloor + 2 &\text{if } n \text{ is odd and }u_n = 1 \\
	\lfloor \frac{n}{2} \rfloor + 1&\text{if } n \text{ is odd and }u_n = 0
	\end{cases}$$
	Thus for all $n \geq 1$ it holds that $P_u(n) \geq P_v(n)$, implying that $11\textswab{t} \in \EL$.

	For minimality, note that $1{\textswab{t}}$ is not prefix normal, since $11$ is a factor of {\textswab{t}}. \qed
\end{proof}

\subsection{Prefix normal forms of Sturmian words.}

Let $w$ be a Sturmian word. As we saw in Sec.~\ref{sec:sturmian}, the only $1$-prefix normal word in the class of Sturmian words with the same slope $\alpha$ is the upper mechanical word $s'_{\alpha,0}=1c_{\alpha}$.

\begin{theorem}\label{thm:sturmian_pnf}
Let $w$ be an irrational mechanical word with slope $\alpha$, i.e.\ a Sturmian word. Then
$\PNF_1(w) = 1c_{\alpha} \text{ and } \PNF_0(w) = 0c_{\alpha},$
where $c_{\alpha}$ is the characteristic word of slope $\alpha$.
\end{theorem}

\begin{proof}
Since the characteristic word $c_{\alpha}$ has the same slope as $w$, we have $\Fac(w)=\Fac(c_{\alpha})$ by Fact~\ref{fact:sturmian-properties}. The abelian complexity of $w$ is constant $2$~\cite{RSZ11}, thus a factor of length $k$ can have either $F^1_w(k)$ or $F^1_w(k)-1$ $1$s. Let us call a factor $u$ of $w$ {\em heavy} if $|u|_1 = F^1_w(k)$, and {\em light} otherwise. We have to show that every prefix of $1c_{\alpha}$ is heavy; this will imply that $1c_{\alpha}$ is the prefix normal form of $w$. It is known~\cite{Lothaire3} that the prefixes of the characteristic word are precisely the reverses of its right special factors, where a factor $u$ is called right special if both $u0$ and $u1$ are factors. Thus, every prefix $v$ of $1c_{\alpha}$ has the form $v = 1u^{\rm rev}$, where both $u1$ and $u0$ are factors of $w$, implying that $|v|_1 = |1u^{\rm rev}|_1 = |u1|_1 = F^1_w(|u|+1),$ therefore $v=1u^{\rm rev}$ is heavy. The fact that $\PNF_0(w) = 0c_{\alpha}$ follows analogously.
\qed \end{proof}

\subsection{Prefix normal forms of binary uniform morphisms}

In~\cite{BSSW16} the authors provide an algorithm which computes the abelian complexity of a morphic word that is the fix point of a binary {\em uniform} morphism, i.e., 
a morphism $\mu$ satisfying $|\mu(0)| = |\mu(1)|$. We refer the reader to~\cite{BSSW16} for the details on this algorithm. In particular, the following theorem is proved in 
\cite{BSSW16}: 

\begin{theorem}[\cite{BSSW16}]
Let $w$ be the fix point of a binary uniform morphism $\mu$. Then, for each $n$ the values $\psi_w(1), \psi_w(2),\dots, \psi_w(n),$ can be computed in $O(n)$ time.
\end{theorem}
 
As an intermediate step in the computation of each $\psi_w(i)$, the algorithm in~\cite{BSSW16} provides the minimum number of $0$s (equivalently, the maximum number of $1$s) in every $i$-length factor of $w$. Obviously the same procedure can be used to obtain  the minimum number of $1$s (equivalently, the maximum number of $0$s) in every $i$-length factor of $w$. 
Therefore, we have the following corollary to the result of \cite{BSSW16}: 

\begin{corollary}
Let $w$ be the fix point of a binary uniform morphism $\mu$. For each $n,$ the prefix of length $n$ of $\PNF_1(w)$ and of $\PNF_0(w)$ can be computed in $O(n)$ time.  
\end{corollary}

\section{Prefix normal words and lexicographic order}\label{sec:lex}

In this section, we study the relationship between lexicographic order and prefix normality. Note that for coherence with the rest of the paper, in the definition of Lyndon words, necklaces, and prenecklaces, we use lexicographically {\em greater} rather than {\em smaller}. Clearly, this is equivalent to the usual definitions up to renaming of characters.

Thus a finite {\em Lyndon word} is one which is lexicographically strictly greater than all of its conjugates: $w$ is Lyndon if and only if for all non-empty $u,v$ s.t.\ $w = uv$, we have  $w >_{\rm lex} vu$. A {\em necklace} is a word which is greater than or equal to all its conjugates, and a {\em prenecklace} is one which can be extended to become a necklace, i.e.\ which is the prefix of some necklace~\cite{Lothaire3,RSW92}. As we saw in the introduction, in the finite case, prefix normality and Lyndon property are orthogonal concepts. However, the set of finite prefix normal words is included in the set of prenecklaces~\cite{BurcsiFLRS17}.

An infinite word is {\em Lyndon} if an infinite number of its prefixes is Lyndon~\cite{SMDS94}. In the infinite case, we have a similar situation as in the finite case. There are words which are both Lyndon and prefix normal: $10^{\omega}, 110(10)^{\omega}$; Lyndon but not prefix normal: $11100(110)^{\omega}$; prefix normal but not Lyndon: $(10)^{\omega}$; and neither of the two: $(01)^{\omega}$.

Next we show that a prefix normal word cannot be lexicographically smaller than any of its suffixes. Let ${\it shift}_i(w) = w_{i}w_{i+1}w_{i+2}\cdots$ denote the infinite word $v$ s.t.\ $w=w_1\cdots w_{i-1}v$, i.e.\ $v$ is the suffix of $w$ starting at position $i$.

\begin{lemma}\label{lemma:pnlex}
Let $w\in \EL_{\rm inf}$. Then $w \geq_{\rm lex} {\it shift}_i(w)$ for all $i\geq 1$.
\end{lemma}

\begin{proof} Assume that there exists a suffix $v = {\it shift}_i(w)$ of $w$ s.t.\ $v>_{\rm lex} w$. Then there is an index $j$ with $v_1\cdots v_{j-1} = w_1\cdots w_{j-1}$ and $v_j > w_j$, implying $v_j=1$ and $w_j = 0$. But then $|w_{i}\cdots w_{i+j-1}|_1 = |v_1\cdots v_j|_1 > |w_1\cdots w_j|_1$, in contradiction to $w\in \EL_{\rm inf}$. \qed
\end{proof}

In the finite case, it is easy to see that a word $w$ is a prenecklace if and only if $w \geq_{\rm lex} v$ for every suffix $v$ of $w$. This motivates our definition of infinite prenecklaces. The situation is the same as in the finite case: prefix normal words form a proper subset of prenecklaces.

\begin{definition}
Let $w \in \{0,1\}^{\omega}$. Then $w$ is an {\em infinite prenecklace} if for all $i\geq 1$, $w \geq_{\rm lex} {\it shift}_i(w)$. We denote by ${\cal P}_{\inf}$ the set of infinite prenecklaces.
\end{definition}

\begin{proposition}
$\EL_{\inf} \subsetneq {\cal P}_{\inf}$.
\end{proposition}

\begin{proof}
The inclusion follows from Lemma~\ref{lemma:pnlex}. An example of a word which is an infinite prenecklace but not prefix normal is $11100(110)^{\omega}$. \qed
\end{proof}

There is another interesting relationship between lexicographic order and the prefix normal forms of an infinite word. In~\cite{Pirillo05}, two words were associated to an infinite binary word $w$, called $\max(w)$ (resp.\ $\min(w)$), defined as the word whose prefix of length $n$ is the lexicographically greatest (resp.\ smallest) $n$-length factor of $w$. It is easy to see that these words always exist. The following was shown in~\cite{Pirillo05}:%
\footnote{The terminology in~\cite{Pirillo05} differs from ours (we are following~\cite{Lothaire3}). In order to help the reader, here we highlight the differences: $(i)$ a periodic Sturmian in~\cite{Pirillo05} is a rational mechanical word, $(ii)$ a proper Sturmian word in~\cite{Pirillo05} is an irrational mechanical word (i.e., a Sturmian word), and $(iii)$ a standard Sturmian word in~\cite{Pirillo05} is a mechanical word with intercept $\tau = \alpha$ (the slope), thus a proper standard Sturmian word is a characteristic Sturmian word $c_\alpha$. Note that all mechanical words in~\cite{Pirillo05} are defined for $n \geq 1$, since the definition of mechanical word is: the lower mechanical word is defined as $ s_{\alpha,\tau} (n) = \lfloor \alpha(n+1) + \tau\rfloor - \lfloor \alpha n  + \tau\rfloor $ for $n \geq 1$, and analogously for the upper mechanical word. Therefore, an intercept $\tau=0$ in~\cite{Pirillo05} is equivalent to an intercept of $\tau=\alpha$ (the slope) in~\cite{Lothaire3}.}

\begin{theorem}[\hspace{1sp}\cite{Pirillo05}]\label{thm:pirillo}
Let $w$ be an infinite binary word. Then
\begin{enumerate}
\item $w$ is (rational or irrational) mechanical with its intercept equal to its slope if and only if $0w \leq_{\rm lex} \min(w) \leq_{\rm lex} \max(w) \leq_{\rm lex} 1w$, and 
\item $w$ is characteristic Sturmian if and only if $\min(w) = 0w$ and $\max(w) = 1w$.
\end{enumerate}

\end{theorem}

\begin{lemma}
Let $w\in\{0,1\}^{\omega}$. Then $\PNF_1(w) \geq_{\rm lex} \max(w)$ and $\PNF_0(w) \leq_{\rm lex} \min(w)$.
\end{lemma}

\begin{proof}
Assume otherwise, and let $w' = \PNF_1(w), v = \max(w)$. If $w'<v$, then there is an index $j$ s.t.\ $w'_1\cdots w'_{j-1} = v_1\cdots v_{j-1}$ and $w'_j=0$ and $v_j = 1$. This implies that $v_1\cdots v_j$ has one more $1$s than $w'_1\cdots w'_j$. But $|w'_1\cdots w'_j|_1 = F^1_w(j)$, a contradiction, since $v_1\cdots v_j$ is a factor of $w$. The second claim follows analogously. \qed
\end{proof}

Finally, from Theorems~\ref{thm:sturmian_pnf} and~\ref{thm:pirillo}, we get the following corollary:

\begin{corollary}\label{coro:pirillo}
Let $w$ be an infinite binary word. Then $w$ is characteristic Sturmian if and only if
$0w = \PNF_0(w) = \min(w) \text{ and } 1w = \PNF_1(w) = \max(w).$
\end{corollary}

\section{On the periodicity and aperiodicity of prefix normal words with respect to minimum density}\label{sec:minimum_density}

In this section, we derive conditions for the periodicity and aperiodicity of prefix normal words with respect to their minimum density.
The following result shows that every ultimately periodic infinite prefix normal word has rational minimum density.

\begin{lemma}
	Let $v$ be an infinite ultimately periodic binary word with minimum density $\delta(v) = \alpha$. Then $\alpha \in {\mathbb Q}$.
\end{lemma}

\begin{proof}
	Let us write $v = ux^{\omega}$ with $x$ not a suffix of $u$.

	For $i=0, 1, \dots , |x|-1$, let $y_i$ be the prefix of length $|u|+i$ of $v$, i.e., $y_i = u x_1 x_2 \cdots x_i.$ Trivially, if for some $i$ we have that
	$\delta(y_i) \leq \delta(v)$ the claim directly follows from $y_i$ being a finite prefix of $v$.

	Let us now assume that for each $i=0,1,\dots |x|-1$ it holds that $\delta(v) < \delta(y_i)$ and let
	$i^*=\min\{ i \mid \delta(y_i) \leq \delta(y_j) \mbox{ for each }j \neq i\}$, hence $\delta(v) < \delta(y_{i^*}).$

	For every $n \geq |u| + |x|$ let $i_n = |u| + ((n - |u|) \bmod |x|)$ and $k_n = \lfloor (n-|u|) / |x| \rfloor$, i.e., $|u| \leq i_n \leq |u| + |x|-1$ and $n = i_n + k_n |x|.$

	Then, we have that
	\begin{equation} \label{lowerbound}
	D_v(n) = \frac{|y_{i_n}|_1 + k_n |x|_1}{|y_{i_n}| + k_n |x|} \geq \min\{\delta(y_{i_n}), \delta(x)\} \geq \min\{\delta(y_{i^*}), \delta(x)\}.
	\end{equation}
	Moreover, we also have that
	\begin{equation} \label{upperbound}
	\lim_{k \to \infty} D_v(|u| + i^*+ k |x|) = \lim_{k \to \infty} \frac{|y_{i^*}|_1 + k |x|_1}{|y_{i^*}| + k |x|} = \delta(x).
	\end{equation}

	We cannot have  $\delta(x) \geq \delta(y_{i^*}),$  since by
	(\ref{lowerbound}) $\delta(y_{i^*})$
	is a rational lower bound on $D_v(n)$ (for each $n \geq 1$) which is achieved by $D_v(|u| + i^*),$ contradicting
	the standing hypothesis $\delta(v) < \delta(y_{i^*}).$

	Therefore, we must have $\delta(x) < \delta(y_{i^*}),$ and from (\ref{lowerbound}) we have $D_v(n) \geq \delta(x)$ and from
	(\ref{upperbound}) we also have that for each $\epsilon > 0$ there exists $k >0$ such that $D_v(|u| + i^*+k|x|) < \delta(x)+\epsilon.$
	Therefore, $\delta(v) = \inf \{ D_v(n) \mid n \geq 1\} = \delta(x),$ which is a rational number, since $x$ is a finite string. \qed
\end{proof}

We now show that, while periodicity is characterized by rational density, the converse is not true. It turns out that for every $\alpha \in (0, 1)$, both rational and irrational, there exists an aperiodic prefix normal word with minimum density $\alpha$. For irrational $\alpha$, this is an easy corollary from Theorem~\ref{thm:characteristic}: since the Sturmian word $1c_{\alpha}$ is prefix normal, and $D(i) \geq \alpha$ for each $i$, therefore, $\delta(1c_{\alpha}) = \alpha$. The next lemma shows how to construct an aperiodic prefix normal word with minimum density $\alpha$ for both rational and irrational $\alpha$.

\begin{lemma}\label{lemma:real_number_density}
Fix $\alpha \in  (0,1)$, and let $(a_n)_{n \in \mathbb{N}}$ be a strictly decreasing infinite sequence of rational numbers from $(0,1)$
converging to $\alpha$.
For each $i = 1, 2, \dots,$ let  the binary word $v^{(i)}$ be defined by
$$
		v^{(i)} = \begin{cases}
		1^{\lceil 10 a_1\rceil} 0^{10-\lceil 10 a_1 \rceil} & i = 1 \\
		\pref_{\flipext^{\omega}(v^{(i-1)})}(k_i |v^{(i-1)}|)0^{\ell_i} & i > 1
		\end{cases}
		\,\,\, 
		$$
where $\ell_i$ defined by
	$$
		\ell_i = \begin{cases}
		10-\lceil 10 a_1 \rceil & i = 1\\
		\left \lfloor k_i \left( \frac{|v^{(i-1)}|_1 - a_i |v^{(i-1)}|}{a_i} \right) \right \rfloor & i > 1,
		\end{cases}
		$$
and $k_i$ is the smallest integer greater than one such that $\ell_i > \ell_{i-1}.$

Then $v = \lim_{i \to \infty} v^{(i)}$ is an aperiodic infinite prefix normal word such that $\delta(v) = \alpha.$
\end{lemma}

Before proving Lemma~\ref{lemma:real_number_density}, in give an example of the words $v^{(i)}$.

\begin{example}
We show the first three steps for the construction of an infinite aperiodic word with minimum density $\alpha=1/3$ (Lemma~\ref{lemma:real_number_density}), using the infinite sequence of rational numbers $a_i = i/(3i-1)$, which tends to $1/3$ for $i \to \infty$. Hence, for $i = 1$, we have $a_1 = 1/2$, $\ell_1 = 5$, and $v_i = 1^5 0^5$ with minimum density $\delta(v_1) = 1/2$. At the next step, $a_2 = 2/5$, and with the values from the previous iteration we can compute $k_2 = 3$ and $\ell_2 = 7$, hence $v_2 = 1^5 0^5 1^5 0^5 1^5 0^5 0^7$, with $\delta(v_2) = 15/37$. At the third iteration, $a_3 = 3/8$, $k_3 = 3$, and $\ell_3 = 9$, therefore $v_3 = 1^5 0^5 1^5 0^5 1^5 0^{12}1^5 0^5 1^5 0^5 1^5 0^{12}1^5 0^5 1^5 0^5 1^5 0^{12} 0^9$, and the minimum density is $\delta(v_3) = 45/120$.
\end{example}

\begin{proof} {\em (of Lemma~\ref{lemma:real_number_density})} 

We will first prove the following claim, giving a number of properties of the sequence of words $v^{(i)}$, and then use these to prove that $v$ is aperiodic and $\delta(v) = \alpha$. 

\medskip 

\noindent
{\em Claim.} The following properties hold: 
\begin{enumerate}
\item  $\delta(v^{(i)}) \geq a_i$  for each $i \geq 1$;
\item $\iota(v^{(i)}) = |v^{(i)}|$ for each $i\geq 1$;
\item $\delta(v^{(i)}) < \delta(v^{(i-1)})$ for each $i \geq 2$;
\item $|v^{(i)}|_1 > |v^{(i-1)}|_1$ for each $i\geq 2$;
\item  $\delta(v^{(i)}) \leq a_i \left(\frac{k_i |v^{(i-1)}|_1}{k_i |v^{(i-1)}|_1 - a_i} \right)$
for each $i\geq 2$.
\end{enumerate}

{\em Proof of the Claim.} By direct inspection we have that properties $1$ and $2$ hold for $v^{(1)}.$ We now argue by induction.
Fix $i > 1$ and let us assume that properties $1$ and $2$ hold for $v^{(i-1)}$. Then, since
$a_i < a_{i-1}$ we have
$$\frac{|v^{(i-1)}|_1}{a_i} > \frac{|v^{(i-1)}|_1}{a_{i-1}} \geq |v^{(i-1)}|,$$ where the last inequality follows from
property $1$ and $2$. Therefore, $\left( \frac{|v^{(i-1)}|_1 - a_i |v^{(i-1)}|}{a_i}\right)$ $> 0$, hence
there exists $k_i > 1$ such that
$\left\lfloor k_i \left( \frac{|v^{(i-1)}|_1 - a_i |v^{(i-1)}|}{a_i}\right) \right\rfloor  > \ell_{i-1}.$ In particular, $\ell_i$ is well defined.

By property 2, we have $\iota(v^{(i-1)}) = |v^{(i-1)}|$ hence by Proposition \ref{flipext-properties}, we have
$D_{\flipext^{\omega}(v^{(i-1)})}(k |v^{(i-1)}|) = \delta(v^{(i-1)})$ and also
$\delta(\pref_{\flipext^{\omega}(v^{(i-1)})}(k_i |v^{(i-1)}|)) = \delta(v^{(i-1)}).$

Moreover, since $\ell_i > 0$ it is not hard to see from the definition of $v^{(i)}$ that
\begin{equation} \label{lemma3:deltavi}
\delta(v^{(i)}) = D_{v^{(i)}}(|v^{(i)}|) = \frac{k_i |v^{(i-1)}|_1}{k_i |v^{(i-1)}| + \ell_i} < \delta(v^{(i-1)}),
\end{equation}
which shows that property 3 and property 2 hold for $v^{(i)}$.
In addition, because of $k_i > 1$ and (by Proposition \ref{flipext-properties}), 
$|v^{(i)}|_1 = |\textstyle{\pref_{\flipext^{\omega}(v^{(i-1)})}(k_i |v^{(i-1)}|)}|_1$ $= k_1 |v^{(i-1)}|_1$, it follows that
property 4 also holds for $v^{(i)}.$ 

The definition of $\ell_i$, together with the well known property $x - 1 < \lfloor x \rfloor \leq x$, imply that
\begin{equation} \label{lemma3:li}
\frac{k_i}{a_i} \left(|v^{(i-1)}|_1 - a_i |v^{(i-1)}| \right)  - 1 < \ell_i \leq k_i \left(\frac{|v^{(i-1)}|_1}{a_i} -  |v^{(i-1)}| \right).
\end{equation}
Using the  right inequality of (\ref{lemma3:li}) in (\ref{lemma3:deltavi}), we have
$\delta(v^{(i)}) \geq a_i$, showing that property 1 holds for $v^{(i)}.$

In addition, using the left inequality of (\ref{lemma3:li}) in (\ref{lemma3:deltavi}), we have
$$\delta(v^{(i)}) \leq a_i \left(\frac{k_i |v^{(i-1)}|_1}{k_i |v^{(i-1)}|_1 - a_i} \right)$$ showing that property 5 holds for
$v^{(i)}.$ The proof of the claim is complete.

\bigskip

In order to see that $v$ is aperiodic, it is enough to observe that $v \neq 0^{\omega}$ and
for each $i \geq 1$ it contains a distinct run of $\ell_i$ $0$s, with $\ell_i$ being a strictly increasing sequence.

To show that $\delta(v) = \alpha,$ we will prove that $\lim_{i \to \infty} \delta(v^{(i)}) = \alpha.$
Since $\lim_{i \to \infty} a_i = \alpha$ and for each $i \geq 1$, $k_i > 1$ and $|v^{(i)}|_1 > |v^{(i-1)}|_1$, we have 
$$\lim_{i \to \infty} a_i \frac{k_i |v^{(i-1)}|_1}{k_i |v^{(i-1)}|_1-a_i} = \lim_{i \to \infty} a_i = \alpha.$$
Hence, from properties $4$ and $5$ of the Claim above, we have the desired result, 
$\lim_{i \to \infty} \delta(v^{(i)}) = \lim_{i \to \infty} a_i = \alpha.$ 

\medskip

This completes the proof of the lemma. \qed

\end{proof}

Summarizing, we have shown the following result.

\begin{theorem}\label{thm:density}
For every $\alpha \in (0,1)$ (rational or irrational) there is an infinite aperiodic prefix normal word of minimum density $\alpha$. On the other hand, for every ultimately periodic infinite prefix normal word $w$, 
the minimum density $\delta(w)$ is a rational number.
\end{theorem}


\section{Conclusion}\label{sec:conclusion}

In this paper, we studied infinite prefix normal words. We gave several results of infinite extensions of finite prefix normal words, and 
we established connections between infinite prefix normal words and other classes of infinite binary words, namely Sturmian words, Lyndon words and max and min words. We provided a complete characterization of prefix normal Sturmian words. Furthermore, we showed that, similar to the finite case, the classes of infinite prefix normal words and Lyndon words are distinct, and that infinite prefix normal words are infinite prenecklaces.

We explored some connections between prefix normal words, prefix normal forms, and abelian complexity. In particular, we showed how to turn balanced and $c$-balanced words without arbitrarily long runs of $1$s into prefix normal words, by prepending a finite number of $1$s. We provided a method to compute the abelian complexity from the prefix normal form of a word, and, for specific cases, we showed how to compute the prefix normal form of a word, given its abelian complexity function. We further applied an existing algorithm to compute the prefix normal form of binary uniform morphisms.

Finally, we gave conditions for the periodicity and the aperiodicity of infinite prefix normal words, according to their minimum density.

\section*{Acknowledgements}
We wish to extend our thanks to the participants of the Workshop on Words and Complexity, which took place in Lyon in February 2018, for exciting discussions and helpful pointers, and to P\'eter Burcsi, who first got us interested in Sturmian words. We also thank the two anonymous reviewers, whose suggestions helped improve the presentation of our results. 
MR is funded by the National Science Foundation (NSF) IIS (Grant No. 1618814), IIBR (Grant No. 2029552) and National Institutes of Health (NIH) R01 (Grant No. HG011392).




\bibliographystyle{abbrv}
\bibliography{PNW}

\begin{thebibliography}{10}

\bibitem{ADKN20}
Peyman Afshani, Ingo van Duijn, Rasmus Killmann, and Jesper~Sindahl Nielsen.
\newblock A lower bound for jumbled indexing.
\newblock In {\em Proceedings of the 2020 {ACM-SIAM} Symposium on Discrete
  Algorithms, ({SODA} 2020)}, pages 592--606, 2020.

\bibitem{AmirCLL_ICAPL14}
Amihood Amir, Timothy~M. Chan, Moshe Lewenstein, and Noa Lewenstein.
\newblock On hardness of jumbled indexing.
\newblock In {\em 41st International Colloquium on Automata, Languages, and
  Programming (ICALP 2014)}, volume 8572 of {\em LNCS}, pages 114--125, 2014.

\bibitem{BG19}
Paul Balister and Stefanie Gerke.
\newblock The asymptotic number of prefix normal words.
\newblock {\em Theoret. Comput. Sci.}, 784:75--80, 2019.

\bibitem{BFR14}
Francine Blanchet-Sadri, Nathan Fox, and Narad Rampersad.
\newblock On the asymptotic abelian complexity of morphic words.
\newblock {\em Advances in Applied Mathematics}, 61:46--84, 2014.

\bibitem{BSSW16}
Francine Blanchet{-}Sadri, Daniel Seita, and David Wise.
\newblock Computing abelian complexity of binary uniform morphic words.
\newblock {\em Theor. Comput. Sci.}, 640:41--51, 2016.
\newblock \href {https://doi.org/10.1016/j.tcs.2016.05.046}
  {\path{doi:10.1016/j.tcs.2016.05.046}}.

\bibitem{BM18}
Alexandre {Blondin Mass{\'{e}}}, Julien de~Carufel, Alain Goupil, M{\'{e}}lodie
  Lapointe, {\'{E}}mile Nadeau, and {\'{E}}lise Vandomme.
\newblock Leaf realization problem, caterpillar graphs and prefix normal words.
\newblock {\em Theoret. Comput. Sci.}, 732:1--13, 2018.

\bibitem{IJFCS12}
P{\'e}ter Burcsi, Ferdinando Cicalese, Gabriele Fici, and {\relax Zs}uzsanna
  Lipt{\'a}k.
\newblock Algorithms for {J}umbled {P}attern {M}atching in {S}trings.
\newblock {\em International Journal of Foundations of Computer Science},
  23:357--374, 2012.

\bibitem{BurcsiCFL12_ToCS}
Peter Burcsi, Ferdinando Cicalese, Gabriele Fici, and {\relax Zs}uzsanna
  Lipt{\'{a}}k.
\newblock On approximate jumbled pattern matching in strings.
\newblock {\em Theory Comput. Syst.}, 50(1):35--51, 2012.

\bibitem{BFLRS20}
P{\'{e}}ter Burcsi, Gabriele Fici, {\relax Zs}uzsanna Lipt{\'{a}}k, Rajeev
  Raman, and Joe Sawada.
\newblock Generating a {Gray} code for prefix normal words in amortized
  polylogarithmic time per word.
\newblock {\em Theor. Comput. Sci.}, 842:86--99, 2020.

\bibitem{BFLRS_CPM14}
P\'eter Burcsi, Gabriele Fici, {\relax Zs}uzsanna Lipt\'ak, Frank Ruskey, and
  Joe Sawada.
\newblock On combinatorial generation of prefix normal words.
\newblock In {\em Proc.\ of the 25th Ann. Symp. on Comb. Pattern Matching (CPM
  2014)}, volume 8486 of {\em LNCS}, pages 60--69, 2014.

\bibitem{BurcsiFLRS17}
P\'eter Burcsi, Gabriele Fici, {\relax Zs}uzsanna Lipt{\'{a}}k, Frank Ruskey,
  and Joe Sawada.
\newblock On prefix normal words and prefix normal forms.
\newblock {\em Theoret. Comput. Sci.}, 659:1--13, 2017.

\bibitem{CK16}
Julien Cassaigne and Idrissa Kabor{\'e}.
\newblock Abelian complexity and frequencies of letters in infinite words.
\newblock {\em Int. Journal of Foundations of Computer Science},
  27(05):631--649, 2016.

\bibitem{ChanL15}
Timothy~M. Chan and Moshe Lewenstein.
\newblock Clustered integer {3SUM} via additive combinatorics.
\newblock In {\em Proc.\ of the 47th Ann. {ACM} on Symp. on Theory of Computing
  ({STOC} 2015)}, pages 31--40, 2015.

\bibitem{CLR18}
Ferdinando Cicalese, {\relax Zs}uzsanna Lipt{\'{a}}k, and Massimiliano Rossi.
\newblock Bubble-flip - {A} new generation algorithm for prefix normal words.
\newblock {\em Theoret. Comput. Sci.}, 743:38--52, 2018.

\bibitem{CLR19}
Ferdinando Cicalese, {\relax Zs}uzsanna Lipt{\'{a}}k, and Massimiliano Rossi.
\newblock On infinite prefix normal words.
\newblock In {\em Proc.\ of the 45th International Conference on Current Trends
  in Theory and Practice of Computer Science ({SOFSEM} 2019)}, pages 122--135,
  2019.

\bibitem{CunhaDGWKS17}
Lu{\'{\i}}s Felipe~I. Cunha, Simone Dantas, Travis Gagie, Roland Wittler, Luis
  Antonio~Brasil Kowada, and Jens Stoye.
\newblock Faster jumbled indexing for binary {RLE} strings.
\newblock In {\em 28th Annual Symposium on Combinatorial Pattern Matching
  ({CPM} 2017)}, pages 19:1--19:9, 2017.

\bibitem{DavisKnuth70}
C.~Davis and D.E. Knuth.
\newblock Number representations and dragon curves, {I, II}.
\newblock {\em J. Recr. Math.}, 3:133--149 and 161--181, 1970.

\bibitem{FL11}
Gabriele Fici and {\relax Zs}uzsanna Lipt{\'a}k.
\newblock On prefix normal words.
\newblock In {\em Proc.\ of the 15th Intern.\ Conf.\ on Developments in
  Language Theory (DLT 2011)}, volume 6795 of {\em LNCS}, pages 228--238.
  Springer, 2011.

\bibitem{FKNP20}
Pamela Fleischmann, Dirk Nowotka, Mitja Kulczynski, and Danny~B{\o}gsted
  Poulsen.
\newblock On collapsing prefix normal words.
\newblock In {\em Proc.\ of the 14th International Conference Language and
  Automata Theory and Applications {(LATA 2020)}}, volume 12038 of {\em LNCS},
  pages 412--424. Springer, 2020.

\bibitem{GagieHLW15}
Travis Gagie, Danny Hermelin, Gad~M. Landau, and Oren Weimann.
\newblock Binary jumbled pattern matching on trees and tree-like structures.
\newblock {\em Algorithmica}, 73(3):571--588, 2015.

\bibitem{GiaGrab_IPL13}
Emanuele Giaquinta and Szymon Grabowski.
\newblock New algorithms for binary jumbled pattern matching.
\newblock {\em Inf. Process. Lett.}, 113(14--16):538--542, 2013.

\bibitem{KaboreK17}
Idrissa Kabor{\'{e}} and Boucar{\'{e}} Kient{\'{e}}ga.
\newblock Abelian complexity of {Thue-Morse} word over a ternary alphabet.
\newblock In {\em Proc.\ of the 11th Int.\ Conf.\ on Combinatorics on Words
  {WORDS 2017}}, volume 10432 of {\em LNCS}, pages 132--143. Springer, 2017.

\bibitem{Lothaire3}
M.~Lothaire.
\newblock {\em Algebraic Combinatorics on Words}.
\newblock Cambridge Univ. Press, 2002.

\bibitem{MadillR13}
Blake Madill and Narad Rampersad.
\newblock The abelian complexity of the paperfolding word.
\newblock {\em Discrete Mathematics}, 313(7):831--838, 2013.
\newblock \href {https://doi.org/10.1016/j.disc.2013.01.005}
  {\path{doi:10.1016/j.disc.2013.01.005}}.

\bibitem{MoosaR_JDA12}
Tanaeem~M. Moosa and M.~Sohel Rahman.
\newblock Sub-quadratic time and linear space data structures for permutation
  matching in binary strings.
\newblock {\em J. Discr. Alg.}, 10:5--9, 2012.

\bibitem{Pirillo05}
Giuseppe Pirillo.
\newblock Inequalities characterizing standard sturmian and episturmian words.
\newblock {\em Theor. Comput. Sci.}, 341(1-3):276--292, 2005.
\newblock \href {https://doi.org/10.1016/j.tcs.2005.04.008}
  {\path{doi:10.1016/j.tcs.2005.04.008}}.

\bibitem{RSZ11}
Gw{\'{e}}na{\"{e}}l Richomme, Kalle Saari, and Luca~Q. Zamboni.
\newblock Abelian complexity of minimal subshifts.
\newblock {\em J. London Math. Society}, 83(1):79--95, 2011.
\newblock \href {https://doi.org/10.1112/jlms/jdq063}
  {\path{doi:10.1112/jlms/jdq063}}.

\bibitem{RSW92}
Frank Ruskey, Carla Savage, and T.M.Y. Wang.
\newblock Generating necklaces.
\newblock {\em J. Algorithms}, 13(3):414--430, 1992.

\bibitem{RSW12}
Frank Ruskey, Joe Sawada, and Aaron Williams.
\newblock Binary bubble languages and cool-lex order.
\newblock {\em J. Comb. Theory, Ser. A}, 119(1):155--169, 2012.

\bibitem{SW12}
Joe Sawada and Aaron Williams.
\newblock Efficient oracles for generating binary bubble languages.
\newblock {\em Electr. J. Comb.}, 19(1):P42, 2012.

\bibitem{SWW17}
Joe Sawada, Aaron Williams, and Dennis Wong.
\newblock Inside the {Binary} {Reflected} {Gray} {Code}: {Flip-Swap} languages
  in 2-{Gray} code order.
\newblock Unpublished manuscript, 2017.

\bibitem{SMDS94}
Rani Siromoney, Lisa Mathew, V.R. Dare, and K.G. Subramanian.
\newblock Infinite {Lyndon} words.
\newblock {\em Inf.\ Proc.\ Letters}, 50:101--104, 1994.

\bibitem{oeis}
N.~J.~A. Sloane.
\newblock The {O}n-{L}ine {E}ncyclopedia of {I}nteger {S}equences.
\newblock Available electronically at \url{http://oeis.org}.

\bibitem{Turek15}
Ondrej Turek.
\newblock Abelian complexity of the {Tribonacci} word.
\newblock {\em J.\ of Integer Sequences}, 18, 2015.

\end{thebibliography}

\end{document}